\newif\ifarxiv 
\arxivtrue 

\newif\ifshort 
\shortfalse

\newif\ifshortt 
\shorttfalse

\documentclass[11pt,a4paper,reqno]{amsart}

\makeatletter

\def\subsection{\@startsection{subsection}{2}%
	\z@{.5\linespacing\@plus.5\linespacing}{.25\linespacing}%
	{\normalfont\bfseries}}

\def\subsubsection{\@startsection{subsubsection}{3}%
	\z@{.5\linespacing\@plus.5\linespacing}{.25\linespacing}%
	{\normalfont\itshape}}

\makeatletter
\def\@settitle{\begin{center}%
  \baselineskip14\p@\relax
	\normalfont\LARGE\scshape  
  \@title
  \end{center}%
}
\makeatother

\newcommand{\key}[1]{\textbf{Keywords.} #1}

\usepackage{color}
\usepackage[dvipsnames]{xcolor}

\usepackage[
bookmarks=true,         
unicode=false,          
pdftoolbar=true,        
pdfmenubar=true,        
pdffitwindow=false,     
pdfstartview={FitV},    
pdftitle={},    
pdfauthor={},     
pdfsubject={},   
pdfcreator={},   
pdfproducer={}, 
pdfkeywords={}, 
pdfnewwindow=true,      
raiselinks	= true,
colorlinks=true,       
linkcolor=MidnightBlue,          
citecolor=Mahogany,        
filecolor=ForestGreen,      
urlcolor=ForestGreen,           
plainpages=false,
pdfpagelabels
]{hyperref}

\date{\today}

\usepackage{eurosym}
\usepackage{amstext} 
\DeclareRobustCommand{\officialeuro}{%
	\ifmmode\expandafter\text\fi
	{\fontencoding{U}\fontfamily{eurosym}\selectfont e}}

\usepackage[T1]{fontenc}     
\usepackage[utf8]{inputenc}	
\usepackage{graphicx, keyval, trig, url} 

\usepackage{epstopdf}

\usepackage{textcomp}
\usepackage{gensymb}
\usepackage[capitalise]{cleveref}
\usepackage{bm}   

\usepackage{blindtext}

\newcommand{\myblue}[1] {{\color{black}#1}}

\usepackage{cite} 

\usepackage{dsfont,amssymb,amsmath, graphicx, amsthm} 
\usepackage{amsfonts,dsfont,mathtools, mathrsfs} 
\usepackage{algorithm,algorithmicx,fancyhdr}
\usepackage[noend]{algpseudocode}

\usepackage{calrsfs}
\DeclareMathAlphabet{\pazocal}{OMS}{zplm}{m}{n}
\let\mathcal\undefined
\DeclareMathAlphabet{\mathcal}{OMS}{cmsy}{m}{n}

\newtheorem{assumption}{Assumption}
\newtheorem{theorem}{Theorem}

\newtheorem{proposition}{Proposition}
\newtheorem{definition}{Definition}

\newtheorem{remark}{Remark}
\newtheorem{corollary}{Corollary}

\newcommand{\Pb}{\mathbb{P}}
\newcommand{\Eb}{\mathbb{E}}

\newcommand{\R}{\mathds{R}}
\newcommand{\N}{\mathds{N}}

\newcommand{\cl}[1]{{\mathcal #1}}

\newcommand{\frk}[1]{{\mathfrak #1}}
\newcommand{\bs}[1]{{\boldsymbol #1}}

\newcommand{\dif}{\mathrm{d}}

\newcommand*{\QEDB}{\hfill\ensuremath{\square}}%

\title[Distributed SMPC for Large-Scale Uncertain Linear Systems]
{Distributed Stochastic Model Predictive Control for
Large-Scale Linear Systems with Private and
Common Uncertainty Sources} 

\author[V.\ Rostampour]{Vahab Rostampour}
\author[T.\ Keviczky]{Tam\'{a}s Keviczky}

\thanks{This research was supported by the Uncertainty Reduction in Smart Energy Systems (URSES) research program funded by the Dutch organization for scientific research (NWO) and Shell under the project Aquifer Thermal Energy Storage Smart Grids (ATES-SG) with grant number 408-13-030.\\
The authors are with Delft Center for Systems and Control, Delft University of Technology, Mekelweg 2, 2628 CD, Delft, The Netherlands. {\tt \{v.rostampour, t.keviczky\}@tudelft.nl}}%

\date{January 4, 2019}

\begin{document}

\begin{abstract}
	This paper presents a distributed stochastic model predictive control (SMPC) approach for large-scale linear systems with private and common uncertainties in a plug-and-play framework. 
	\myblue{
		Using the so-called scenario approach, the centralized SMPC involves formulating a large-scale finite-horizon scenario optimization problem at each sampling time, which is in general computationally demanding, due to the large number of required scenarios. 
		We present two novel ideas in this paper to address this issue. 
		We first develop a technique to decompose the large-scale scenario program into distributed scenario programs that exchange a certain number of scenarios with each other in order to compute local decisions using the alternating direction method of multipliers (ADMM).
	}
	We show the exactness of the decomposition with a-priori probabilistic guarantees for the desired level of constraint fulfillment for both uncertainty sources. 
	As our second contribution, we develop an inter-agent soft communication scheme based on a set parametrization technique together with the notion of probabilistically reliable sets to reduce the required communication between the subproblems. 
	We show how to incorporate the probabilistic reliability notion into existing results and provide new guarantees for the desired level of constraint violations. 
	\myblue{
		Two different simulation studies of two types of systems interactions, dynamically coupled and coupling constraints, are presented to illustrate the advantages of the proposed framework.
	}
\end{abstract}

\maketitle

\key{
	Stochastic MPC, 
	Scenario MPC, 
	Distributed MPC, 
	Distributed Stochastic MPC,
	Distributed Scenario MPC,
	Decentralized Scenario MPC,
	Plug-and-Play Framework.
}


\section{Introduction}

Stochastic
model predictive control (SMPC) has attracted significant attention in the recent control literature, due to its ability to provide an alternative, often less conservative way to handle uncertain systems. 
SMPC takes into account the stochastic characteristics of the uncertainties and thereby the system constraints are treated in a probabilistic sense, i.e., using chance constraints \ifshortt\cite{kouvaritakis2004recent}\else\cite{kouvaritakis2004recent, mesbah2016stochastic}\fi.
SMPC computes an optimal control sequence that minimizes a given objective function subject to the uncertain system dynamics model and chance constraints in a receding horizon fashion \cite{primbs2007soft}.
Chance constraints enable SMPC to offer an alternative approach to achieve a less conservative solution compared to robust model predictive control (MPC) \cite{bemporad1999robust}, since it directly incorporates the tradeoff between constraint feasibility and control performance.

\ifshortt

Distributed MPC has been an active research area in the past decades, due to its applicability to handle large-scale dynamic systems with state and input constraints.
\else

Distributed MPC has been an active research area in the past decades, due to its applicability in different domains such as power networks \cite{venkat2008distributed}, chemical plants \cite{rostampour2015stochastic}, process control \cite{ydstie2002new}, and building automation \cite{ma2012dmpc}. 
For such large-scale dynamic systems with state and input constraints, distributed MPC is an attractive control scheme.
\fi
In distributed MPC one replaces large-scale optimization problems stemming from centralized MPC with several smaller-scale problems that can be solved in parallel.
These problems make use of information from other subsystems to formulate finite-horizon optimal control problems.
In the presence of uncertainties, however, the main challenge in the formulation of distributed MPC is how the controllers should exchange information through a communication scheme among subsystems (see, e.g., \cite{lavaei2008control}, and references therein).
This highlights the necessity of developing distributed control strategies to cope with the uncertainties in subsystems while at the same time minimizing information exchange through a communication framework.

\ifarxiv
\subsection*{Related Works}
\else
\fi

In order to handle uncertainties in distributed MPC, some approaches are based on robust MPC \ifshortt\cite{richards2007robust}\else\cite{richards2007robust, conte2013robust}\fi. 
Assuming that the  uncertainty is bounded, a robust optimization problem is solved at each sampling time, leading to a control law that satisfies the constraints for all  admissible values of the uncertainty.
The resulting solution using such an approach tends to be conservative in many cases.
Tube-based MPC, see for example \cite{cannon2011stochastic} and the references therein, was considered in a plug-and-play decentralized setup in \cite{riverso2013plug}, and it has been recently extended to distributed control systems \cite{dai2016distributed} for a collection of linear stochastic subsystems with independent dynamics.
\ifshortt
\else
While in \cite{dai2016distributed} coupled chance constraints were considered separately at each sampling time, in this paper we consider a chance constraint on the feasibility of trajectories of dynamically coupled subsystems.
Our approach is motivated by \cite{riverso2013plug} to reduce the conservativeness of the control design.
\fi
\ifshortt 
\else
Other representative approaches for SMPC of a single stochastic system include affine parametrization of the control policy \cite{hokayem2012stochastic}, the randomized (scenario) approach \cite{prandini2012randomized, calafiore2013stochastic, schildbach2015scenario, lorenzen2015scenario}, and the combined randomized and robust approach \cite{margellos2014road, zhang2013stochastic, rostampour2016robust}.
None of these approaches, to the best of our knowledge, have been considered in a distributed control setting.
\fi

This paper aims to develop a systematic approach to distributed SMPC using the scenario MPC technique.
Scenario MPC approximates SMPC via the so-called scenario (sample) approach \ifshortt\cite{campi2008exact}\else\cite{calafiore2006scenario, campi2008exact}\fi, and if the underlying optimization problem is convex with respect to the decision variables, finite sample guarantees can be provided.
Following such an approach, the computation time for a realistic large-scale system of interest becomes prohibitive, due to the fact that the number of samples to be extracted tends to be high, and consequently leads to a large number of constraints in the resulting optimization problem.
\ifshortt
\myblue{
	To summarize the main contributions of this paper are as follows:
	
	(a) By reformulating a centralized scenario optimization problem into a decomposable scenario program, we first present the exactness of decomposition and then, quantify the level of robustness of resulting solutions by providing new a-priori probabilistic guarantees for the desired level of constraint fulfillment under some mild conditions.  
	
	(b) To solve the proposed decomposable scenario program, we develop a so-called distributed scenario exchange scheme between the neighboring agents using the ADMM. Such a new scheme is used to provide a distributed scenario MPC framework to handle both dynamically coupled and coupling constraints network of agents. 
	
	(c) To reduce the communication in the proposed distributed SMPC, we develop a novel inter-agent soft communication scheme based on a set parametrization technique together with the notion of probabilistically reliable set to reduce the required communication between the subproblems. We show how to incorporate the probabilistic reliability notion into existing results and provide new guarantees for the desired level of constraint violations.
	
	(d) Using the so-called soft communication scheme, we present a distributed SMPC such that the neighboring agents just need to communication a set together with the level of reliability of the set. Each agent then solves a \textit{local robust-communication scenario program} at each sampling time. We also establish a practical plug-and-play (PnP) distributed SMPC framework that considers network changes by agents which want to join or leave the network.
}

\else
To  overcome  the  computational  burden caused by the large number of constraints, in \cite{long2014scenario, liu2016scenario} a heuristic sample-based approach was used in an iterative distributed fashion via dual decomposition such that all subsystems collaboratively optimize a global performance index.   
In another interesting work \cite{calafiore2013random}, a multi-agent consensus algorithm was presented  to achieve  consensus  on  a  common  value  of the decision vector subject to random constraints such that a probabilistic bound on the tails of the consensus violation was also established.
However, in most of the aforementioned references the aim to reduce communication among subsystems, which we refer to as agents, has not been addressed.

\ifarxiv
\subsection*{Contributions}
\else
\fi

\myblue{
	Our work in this paper differs from the above references in two important aspects which have not been, to the best of our knowledge, considered in literature.
	We develop two different distributed scenario MPC frameworks.
	In the first version, each agent generates its own scenarios of the local uncertainty sources, and takes into consideration an estimation of the neighboring possible state trajectories, such that using the so-called alternating direction method of multipliers (ADMM), all the agents start to improve their local estimation until they agree on the consistency of exchange scenarios in the proposed distributed setting. 
	In the second version, we develop a new set-based communication setup together with the notion of probabilistically reliable information to reduce the communications between agents required by the proposed first version of distributed scenario MPC framework. To quantify the error introduced by such a new communication scheme, we incorporate the probabilistic reliability notion into existing results and provide new guarantees for the desired level of constraint violations.  
	To summarize the main contributions of this paper are as follows:
	
	\begin{enumerate}
		
		\item By reformulating a centralized scenario optimization problem into a decomposable scenario program, we first present the exactness of decomposition and then, quantify the level of robustness of resulting solutions by providing new a-priori probabilistic guarantees for the desired level of constraint fulfillment under some mild conditions.  
		
		\item To solve the proposed decomposable scenario program, we develop a so-called distributed scenario exchange scheme between the neighboring agents using the ADMM. Such a new scheme is used to provide a distributed scenario MPC framework to handle both dynamically coupled and coupling constraints network of agents. 
		
		\item To reduce the communication in the proposed distributed SMPC, we develop a novel inter-agent soft communication scheme based on a set parametrization technique together with the notion of probabilistically reliable set to reduce the required communication between the subproblems. We show how to incorporate the probabilistic reliability notion into existing results and provide new guarantees for the desired level of constraint violations.
		
		\item Using the so-called soft communication scheme, we present a distributed SMPC such that the neighboring agents just need to communication a set together with the level of reliability of the set. Each agent then solves a \textit{local robust-communication scenario program} at each sampling time. We also establish a practical plug-and-play (PnP) distributed SMPC framework that considers network changes by agents which want to join or leave the network.
		
	\end{enumerate}
	
}

\fi

It is important to highlight that two major difficulties arising in SMPC, namely recursive feasibility \cite{lofberg2012oops} and stability, are not in the scope of this paper, and they are subject of our ongoing research work.
Thus, instead of analyzing the closed-loop asymptotic behavior, in this paper we focus on individual SMPC problem instances from the optimization point of view and derive probabilistic guarantees for constraint fulfillment in a distributed setting.
\ifarxiv
\else
\myblue{The reader is also referred to the extended version of this paper available in \cite{rostampour2018distributedArXiV}.}
\fi

\ifshort
\else

\ifarxiv
\vspace{0.25cm}
\subsection*{Structure}
\else
\fi

\myblue{
	The structure of this paper is as follows.
	\Cref{prob_state} describes a mathematical model of the control system dynamics together with formulating a large-scale centralized scenario MPC.
	In \Cref{dist_prob_state}, we first reformulate the centralized problem into a decomposable scenario problem using a decomposition technique.
	We then provide the results on equivalent relations and analyze robustness of the obtained solutions  at each sampling time using the proposed distributed scenario exchange scheme in \Cref{dist_smpc} within the proposed framework for the distributed scenario MPC.
	\Cref{inter_agent} introduces a novel inter-agent communication scheme between the subproblems, namely \textit{soft communications}, and then proceeds to quantify the robustness of the proposed schemes. 
	Using the proposed soft communication scheme, in \Cref{PnP} we establish a practical PnP distributed SMPC framework.
	\Cref{sim} presents two different simulation studies to illustrate the functionality of our theoretical achievements, whereas in \Cref{final}, we conclude this paper with some remarks and future work.
}

\fi

\section*{Notations}

$\R, \R_{+}$ denote the real and positive real numbers, and $\N, \N_{+}$ the natural and positive natural numbers, respectively.
We operate within the $n$-dimensional space $\R^n$ composed of column vectors $u, v \in \R^n$.
The Cartesian product over $n$ sets $\cl{X}_1, \cdots, \cl{X}_n$ is given by: $\prod_{i=1}^n \cl{X}_i = \cl{X}_1 \times \cdots \times \cl{X}_n = \{(x_1, \cdots, x_n) \, : \, x_i \in \cl{X}_i\}$.
The cardinality of a set $\cl{A}$ is denoted by $|\cl{A}| = A$.
We denote a block-diagonal matrix with blocks $X_i$, $i\in\{1,\cdots, n\}$, by $\text{diag}_{i\in\{1,\cdots, n\}}(X_i)$, and a  vector consisting of stacked sub vectors $x_i$, $i\in\{1,\cdots, n\}$, by $\text{col}_{i\in\{1,\cdots, n\}}(x_i)$.

Given a metric space $\Delta$, its Borel $\sigma$-algebra is denoted by $\frk{B}(\Delta)$.
Throughout the paper, measurability always refers to Borel measurability.
In a probability space $(\Delta, \frk{B}(\Delta),\Pb)$, we denote the $N$-Cartesian product set of $\Delta$ by $\Delta^N$ and the respective product measure by $\Pb^N$.

\section{Problem Formulation}\label{prob_state}

\ifshortt
\else

This section provides an overview of the control problem statement.
We first describe the dynamics of a large-scale uncertain linear system together with input and state constraint sets and the control objective.
We then formulate a centralized SMPC for such a large-scale control system problem.
Finally, a tractable reformulation based on the scenario MPC \cite{schildbach2014scenario} together with theoretical connections are provided.

\fi

Consider a discrete-time uncertain linear system with additive disturbance in a compact form as follows:
\begin{align}\label{dyn_LS}
{x}_{k+1} = A(\delta_k) {x}_{k} + B(\delta_k) {u}_{k} + C(\delta_k) w_{k} \ ,
\end{align}
with a given initial condition ${{x}}_0\in\R^n$.
Here $k\in\cl{T} := \{0, 1, \cdots, T-1\}$ denotes the time instance,   
${x}_k\in\cl{X}\subset \R^n$ and ${u}_k\in\cl{U}\subset \R^m$  correspond to the state and control input, respectively, and $w_k\in \R^p$ represents an additive disturbance.
The system matrices $A(\delta_k) \in \R^{n\times n}$ and $B(\delta_k) \in \R^{n\times m}$ as well as $C(\delta_k) \in \R^{n\times p}$ are random, since they are known functions of an uncertain variable $\delta_k$ that influences the system parameters at each time step $k$.
\myblue{It is important to mention that given the initial measured state variable $x_0$, the predicted future states are modeled using \eqref{dyn_LS} at each sampling time $k$. 
	For avoiding crowded notation, when there is no ambiguity, in the rest of the paper we will drop the conventional index for the predicted time steps, i.e., $k+\ell|k$, and simply use $k+\ell$.
	We also define $\cl{T}_+ := \{1, \cdots, T\}$ to denote the future predicted time steps. 
}

\begin{assumption}
	Random variables $\bs{w} := \{w_k\}_{k\in\cl{T}}$ and $\bs{\delta} := \{\delta_k\}_{k\in\cl{T}}$ are defined on probability spaces $(\cl{W},\frk{B}(\cl{W}),\Pb_{\bs{w}})$ and $({\Delta},\frk{B}({\Delta}),\Pb_{\bs{\delta}})$, respectively.
	$\bs{w}$ and  $\bs{\delta}$ are two independent random processes, where 	 $\Pb_{\bs{w}}$ and $\Pb_{\bs{\delta}}$  are two different probability measures defined over $\cl{W}$ and ${\Delta}$, respectively, and $\frk{B}(\cdot)$ denotes a Borel $\sigma$-algebra. 
	The support sets $\cl{W}$ and $\Delta$ of $\bs{w}$ and $\bs{\delta}$, respectively, together with their probability measures $\Pb_{\bs{w}}$ and $\Pb_{\bs{\delta}}$  are entirely generic.  
	In fact, $\cl{W}$, $\Delta$ and $\Pb_{\bs{w}}$, $\Pb_{\bs{\delta}}$  do not need to be known explicitly.  
	Instead, the only requirement is availability of a "sufficient number" of samples, which will become concrete in later parts of the paper. 
	Such samples can be for instance obtained by a model learned from available historical data \cite{papaefthymiou2008mcmc}.
\end{assumption}

The system in \eqref{dyn_LS} is subject to constraints on the system state trajectories and control input.
Consider the state and control input constraint sets to be compact convex in the following form
\begin{equation}\label{feasible_set_LS}
\begin{aligned}
\mathcal{X} := \{{x} \in \mathbb{R}^{n} : \, G \, {x} \leq g\}, 
\mathcal{U} := \{{u} \in \mathbb{R}^{m} : \, H \, {u} \leq h\},
\end{aligned}
\end{equation}
where $G \in \mathbb{R}^{q \times n}$, $g \in \mathbb{R}^{q}$, and $H \in \mathbb{R}^{r \times m}$, $h \in \mathbb{R}^{r}$.
Keeping the state inside a feasible set $\mathcal{X}\subset \R^n$ for the entire prediction horizon may be too conservative and result in loss of performance. 
In particular, this is the case when the
best performance is achieved close to the boundary of $\mathcal{X}$, and
thus, constraint violations will be unavoidable due to the fact
that the system parameters in  \eqref{dyn_LS} are imperfect and uncertain. 
To tackle such a problem, we will consider chance constraints on the state trajectories to avoid violation of the state variable constraints probabilistically even if the disturbance $\bs{w}$ or uncertainty $\bs{\delta}$ has unbounded support.
Notice that a robust problem formulation \cite{bemporad1999robust} cannot cope with problems having an unbounded disturbance set.

In order to find a stabilizing full-information controller that leads to admissible control inputs $\bs{u} := \{{u}_k\}_{k\in\cl{T}}$ and satisfies the state constraints, we follow the traditional MPC approach.
The design relies on the standard assumption of the existence of a suitable pre-stabilizing control law, see, e.g., \cite[Proposition 1]{riverso2013plug}. 
To cope with the state prediction under uncertainty and disturbance, we employ a parametrized feedback policy \cite{hokayem2012stochastic}  and split the control input:
${u}_{k} = K {x}_{k} + {v}_{k}$ with ${v}_{k}\in\R^{m}$ as a free correction input variable to compensate for disturbances.

The control objective is to minimize a cumulative quadratic stage cost of a finite horizon cost $J(\cdot): \R^n\times \R^m \rightarrow \R$  that is defined as follows:
\begin{align}
J(\bs{x}, \bs{u}) = \mathbb{E} \left[ \sum_{k=0}^{T-1} \left({x}_{k}^\top Q {x}_{k}+ {u}_{k}^\top R {u}_{k}\right) + {x}_{T}^\top P {x}_{T} \right] \,,
\end{align}
with $Q \in \mathbb{R}^{n \times n}_{\succeq 0}$, and $R \in \mathbb{R}^{m \times m}_{\succ 0}$.
Consider $\bs{x}:= \{x_k\}_{k\in\cl{T}}$, $(A, Q^{\frac{1}{2}})$ to be detectable and $P$ to be the solution of the discrete-time Lyapunov equation:
\begin{align}\label{feedback_law}
\Eb[\,A_{cl}(\delta_k)^\top P A_{cl}(\delta_k)\,] + Q + K^\top R K - P \preceq 0\,,
\end{align}
for the closed-loop system with $A_{cl}(\delta_k) = A(\delta_k) + B(\delta_k) K$.
Each stage cost term is taken in expectation $\Eb[\cdot]$, since the argument ${x}_{k}$ is a random variable.  
Using $\bs{v} = \{v_k\}_{k\in\cl{T}}$, consider now the following stochastic control problem:
\begin{subequations}
	\label{scp_ccp_ls}
	\begin{align}
	\min_{\bs{v}\in\R^{Tm}} & \quad J(\bs{x}, \bs{u}) \label{ccp_cost}  \\
	\text{s.t.} \ \ & \quad {x}_{k+1} = A(\delta_k) {x}_{k} + B(\delta_k) {u}_{k} + C(\delta_k) {w}_{k} \,,  \\
	& \quad \Pb [\, {x}_{k+\ell} \in \cl{X} \,,\, \ell \in \cl{T}_+\,] \geq 1 - \varepsilon\,,  \label{ccp_ls}\\
	& \quad {u}_k = K {x}_{k} + {v}_{k} \in \cl{U} \ ,\quad \forall k\in\cl{T} \,, 
	\end{align}
\end{subequations}
where ${{x}}_0$ is initialized based on the measured current state, and $\varepsilon\in(0,1)$ is the admissible state constraint violation parameter of the large-scale system \eqref{dyn_LS}.
\ifshortt
\else
The objective function is assumed to be a quadratic function; however, this is not a restriction and any generic convex  function can be chosen instead.
It is important to mention that the parameters of constraint sets, $\cl{X}$, $\cl{U}$, and the objective function $J(\cdot)$ can be time-varying with respect to the sampling time $k\in\cl{T}$. For the clarity of our problem formulation, we assume time-invariance.
\fi 
The state trajectory ${x}_{k+\ell}\,, \forall \ell\in\cl{T}_+$, has a dependency on the random variables $\bs{w}$ and $\bs{\delta}$, and thus, the chance constraint can be interpreted as follows:
the probability of violating the state constraint at the future time step $\ell\in\cl{T}_+$ is restricted to $\varepsilon$, given that the state of the system in \eqref{dyn_LS} is measurable at each time step ${k\in\cl{T}}$.
Even though $\mathcal{U}$ and $\mathcal{X}$ are compact convex sets, due to the chance constraint on the state trajectory, the feasible set of the optimization problem in \eqref{scp_ccp_ls} is a non-convex set, in general.

\ifshortt
\else

\begin{remark}
	Instead of the chance constraint on the state trajectory of form \eqref{ccp_ls}, one can also bound the average rate of state constraint violations \cite{schildbach2014scenario}.
	Moreover, one can also define the cost function \eqref{ccp_cost} as a desired  quantile of the sum of discounted stage costs ("value-at-risk"), instead of the sum of expected values. 
	Instead of a state feedback law, one can also consider a nonlinear disturbance parametrization feedback policy over the prediction horizon, similar to \cite{rostampour2016robust}, using the scenario-based approximation.
	Such a parametrization does not affect the convexity of the resulting optimization \cite{prandini2012randomized}.
\end{remark}

\fi

To handle the chance constraint \eqref{ccp_ls}, we recall a scenario-based approximation \cite{schildbach2014scenario}.
$w_k$ and $\delta_k$ at each sampling time $k\in\cl{T}$ are not  necessarily  independent and identically distributed  (i.i.d.).
In particular,  they  may  have time-varying distributions and/or be correlated in time.	
We assume that a "sufficient number" of i.i.d. samples of the disturbance $\bs{w}\in\cl{W}$ and $\bs{\delta}\in\Delta$ can be obtained either empirically or by a random number generator.
We denote the sets of given finite samples (scenarios) with $\cl{S}_{\bs{w}}:= \{\bs{w}^{1}, \cdots, \bs{w}^{S}\}\in\cl{W}^{S}$ and $\cl{S}_{\bs{\delta}}:= \{\bs{\delta}^{1}, \cdots, \bs{\delta}^{S}\}\in\Delta^{S}$, respectively.
\myblue{
	Following the approach in \cite{rostampour2015stochastic}, we approximate the expected value of the objective function empirically by averaging the value of its argument for some number of different scenarios, which plays a tuning parameter role. 
	Using $\bar{S}$ as the tuning parameter, consider $\bar{S}$ number of different scenarios of $\bs{w}$ and $\bs{\delta}$ to build 
	\[\bar{\cl{S}}_{\bs{w}, \bs{\delta}}  = 
	\left\{ (\bs{w}^{s},\bs{\delta}^{s}) \, : \, \bs{w}^{s} \in \mathcal{W} \, , \,  \bs{\delta}^{s} \in \Delta \ , s = 1, \cdots, \bar{S} \right\}\,,
	\]
	which has the cardinality $|\bar{\cl{S}}_{\bs{w}, \bs{\delta}}| = \bar{S}$. 
	We then approximate the cost function empirically as follows: 
	\begin{align*}
	J(\bs{x}, \bs{u}) = \Eb_{(\bs{w}, \bs{\delta})} \left[ \sum_{k=0}^{T-1} V({x}_{k}({w}_{k}, {\delta}_{k}), {u}_{k})\right] \approx \frac{1}{\bar{S}} \sum_{(\bs{w}^{s}, \bs{\delta}^{s}) \in \bar{\cl{S}}_{\bs{w},\bs{\delta}} } \sum_{k=0}^{T-1} V({x}_{k}({w}_{k}^{s}, {\delta}_{k}^{s}), {u}_{k}) \ ,
	\end{align*}
	where $V({x}_{k}({w}, {\delta}_{k}), {u}_{k})$ represents a compact notation of the objective function
	$\left({x}_{k}({w}_{k}, {\delta}_{k})^\top Q {x}_{k}({w}_{k}, {\delta}_{k})+ {u}_{k}^\top R {u}_{k}\right) + {x}_{T}({w}_{k}, {\delta}_{k})^\top P {x}_{T}({w}_{k}, {\delta}_{k})$.
	Notice that ${x}_{k}({w}_{k}, {\delta}_{k})$ indicates the dependency of the state variables on the random variables.  
}

We are now in a position to formulate an approximated version of the proposed stochastic control problem in \eqref{scp_ccp_ls} using the following finite horizon scenario program:
\begin{subequations}
	\label{scp_ls}
	\begin{align}
	\min_{\bs{v}\in\R^{T{m}}}& \ J(\bs{x}, \bs{u}) \\
	\text{s.t.} \ \ & \  {x}_{k+1}^{s} = A(\delta_{k}^{s}) {x}_{k}^{s} + B(\delta_{k}^{s}) {u}_{k}^{s} + C(\delta_{k}^{s}) {w}_{k}^{s} \,,   \\
	& \  {x}_{k+\ell}^{s} \in \cl{X} \,,\, \ell \in \cl{T}_+\,,\, \forall\bs{w}^{s} \in {\cl{S}}_{\bs{w}} \,,\, \forall\bs{\delta}^{s} \in {\cl{S}}_{\bs{\delta}}\,,  \\
	&	\  {u}_k^{s} = K {x}_k^{s} + {v}_{k} \in \cl{U} \ ,\quad \forall k\in\cl{T} \,, \label{cons_u} 
	\end{align}
\end{subequations}
where superscript $s$ indicates a particular sample realization.
The solution of \eqref{scp_ls} is the optimal control input sequence $\bs{v}^* = \{v_{k}^{*}, \cdots, v_{k+T-1}^{*}\}$.
Based on the MPC paradigm, the current input is implemented as $u_{k} := K x_k + v_{k}^{*}$ and we proceed in a receding horizon fashion. 
This means that the problem \eqref{scp_ls} is solved at each time step $k$ by using the current measurement of the state $x_{k}$. 
Note that new scenarios are needed to be generated at each sampling time $k\in\N_+$.
\myblue{
	It is important to note that the proposed centralized scenario program \eqref{scp_ls} should be feasible and its feasibility domain has to be a nonempty domain. This is a technical requirement (\cite[Assumption 1]{campi2008exact}) to be able to use such a randomization technique. 
	In case of infeasible solution, we need to generate new set of scenarios and resolve the problem.
}

\myblue{
	\begin{assumption}\label{ass_scp_ls_feasible}
		The proposed centralized scenario program \eqref{scp_ls} is feasible and its feasibility domain has a nonempty domain.
	\end{assumption}

	Following Assumption~\ref{ass_scp_ls_feasible} and based on the Weierstrass' theorem \cite[Proposition~A.8]{bertsekas1989parallel}, the proposed centralized scenario program \eqref{scp_ls} admits  at  least  one  optimal solution, and therefore the Slater’s constraint qualification \cite{slater1959lagrange}
	holds for the proposed problem formulations.
}
The key features of the proposed optimization problem \eqref{scp_ls} are as follows:
1) there is no need to know the probability measures $\Pb_{\bs{w}}$ and $\Pb_{\bs{\delta}}$ explicitly, only the capability of obtaining random scenarios is enough,
2) formal results to quantify the level of approximations are available. In particular, the results follow the so-called scenario approach \cite{campi2008exact}, which allows to bound a-priori the violation probability of the obtained solution via \eqref{scp_ls}.

In the following theorem, we restate the explicit theoretical bound of \cite[Theorem 1]{campi2008exact}, which measures the finite scenarios behavior of \eqref{scp_ls}. 
\begin{theorem}\label{scenario_thm}
	Let $\varepsilon\,,\beta \in (0,1)$ and $S \geq \mathsf{N}(\varepsilon, \beta, Tm)$, where
	\begin{equation*}
	\begin{aligned}
	\mathsf{N}(\varepsilon, \beta, Tm) := \min\left\{ \mathbf{N} \in \N \Big| \sum_{i=0}^{Tm-1} \binom{\mathbf{N}}{i} \varepsilon^{i} (1-\varepsilon)^{\mathbf{N}-i} \leq \beta \right\}.
	\end{aligned}
	\end{equation*}
	If the optimizer of problem \eqref{scp_ls} $\bs{v}^*\in\R^{Tm}$ is applied to the discrete-time dynamical system \eqref{dyn_LS} for a finite horizon of length $T$, then, with at least confidence $1-\beta$, the original constraint \eqref{ccp_ls} is satisfied for all $k\in\cl{T}$ with probability more than $1-\varepsilon$.
\end{theorem}
It was shown in \cite{campi2008exact} that the above bound is tight.
The interpretation of \Cref{scenario_thm} is as follows: when applying $\bs{v}^*$ in a finite horizon control problem, the probability of constraint violation of the state trajectory remains below $\varepsilon$ with confidence $1-\beta$: 
\begin{align*}
\Pb^{S}\left[\, \cl{S}_\bs{w}\in\cl{W}^{S}, \cl{S}_\bs{\delta}\in\Delta^{S}\,:\, \text{Vio}(\bs{v}^*) \leq \varepsilon \,\right] \geq 1-\beta \,,
\end{align*}
with $\text{Vio}(\bs{v}^*) := \Pb[\, \bs{w}\in\cl{W}, \bs{\delta}\in\Delta \,:\, {x}_{k+\ell} = A_{cl}(\delta_k) {x}_{k} + B(\delta_k) {v}_{k}^* + C(\delta_k) {w}_{k} \notin\cl{X} \,,\, \ell \in \cl{T}_+ \big\vert \, {x}_k = {{x}}_0 \,] \,,$
where $A_{cl}(\delta_k) = A(\delta_k) + B(\delta_k) K$.
\myblue{
	It is worth to mention that while the proposed constraint on the control input \eqref{cons_u} is also met in a probabilistic sense for all prediction time steps, except at the initial time step, due to the feedback parametrization and the nature of the scenario approach that appears in the proposed optimization problem \eqref{scp_ls}.
	At the initial time step such a constraint \eqref{cons_u} is deterministic, and the superscript $s$ can be dropped, as there is only one measured current state $x_k$. This holds for all the following proposed formulations.
}

\ifshortt
\else

\begin{remark}
	One can obtain an explicit expression for the desired number of scenarios $S$ as in \cite{alamo2015randomized}, where it is shown that given $\varepsilon, \beta \in (0,1)$ and $e$ the Euler constant, then $S\geq \frac{e}{e-1}\frac{1}{\varepsilon} \left(Tm + \ln \frac{1}{\beta}\right)$. 
	It is important to note that $S$ is used to construct the sets of scenarios, $\cl{S}_\bs{w}$, $\cl{S}_\bs{\delta}$ to obtain a probabilistic guarantee for the desired level of feasibility, while the number of scenarios $\bar{S}$ is just a tuning variable to approximate the objective function empirically.
\end{remark}

\fi

\ifshort
\else

We formulated a large-scale SMPC \eqref{scp_ccp_ls} together with a tractable reformulation based on the proposed centralized scenario MPC \eqref{scp_ls}.
Fig.~\ref{fig:smpc} shows a pictorial representation of \eqref{scp_ls} as a large-scale network of interconnected agents to summarize this section.
It is worth mentioning that such a large-scale SMPC \eqref{scp_ccp_ls} is initially proposed in \cite{rostampour2017probabilistic} for a network of interconnected buildings in smart thermal grids.
In the following section, we will provide a distributed framework to solve the proposed problem in \eqref{scp_ls} by decomposing the large-scale system dynamics \eqref{dyn_LS}.

\begin{figure}
	\centering
	\includegraphics[width=0.65\textwidth]{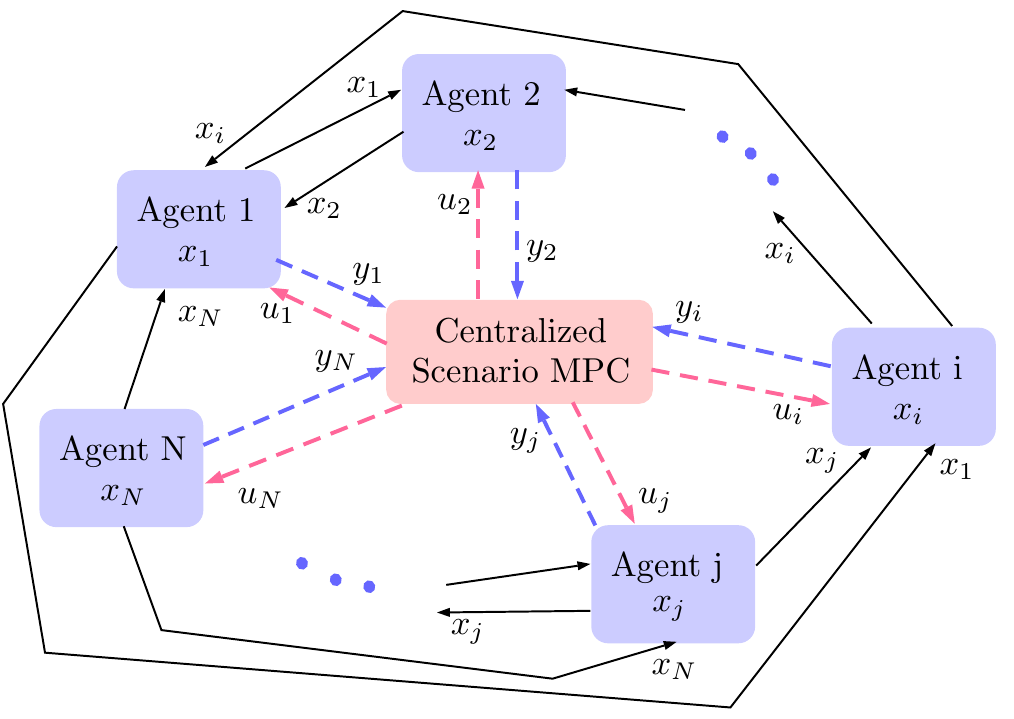}
	\caption{Centralized scenario MPC that corresponds to the problem \eqref{scp_ls}. 
		The measurement variable $y_i$ is the full state information $x_i$ which is sent to the centralized controllers for all agents $i\in\mathcal{N}$.}
	\label{fig:smpc}
\end{figure}

\fi

\section{Decomposed Problem Reformulation}\label{dist_prob_state}

\ifshortt
\else

In this section, we describe a decomposition technique to partition the large-scale system dynamics in \eqref{dyn_LS}.
By taking into consideration two possible uncertainty sources, namely private (local) and common uncertainties, for the resulting network of interconnected subsystems (agents), we provide the theoretical links to the results that it is provided in the previous section.

\fi

Consider a partitioning of the system dynamics \eqref{dyn_LS} 
through a decomposition into $N$ subsystems and let $\cl{N}=\{1, 2, \cdots, N\}$ be the set of subsystem indices.
The state variables $x_k$, control input signals $u_k$ and the additive disturbance $w_k$ can be considered as ${x}_k = \text{col}_{i\in\cl{N}}(x_{i,k})$, ${u}_k = \text{col}_{i\in\cl{N}}(u_{i,k})$ and ${w}_k = \text{col}_{i\in\cl{N}}(w_{i,k})$, respectively, where $x_{i,k}\in\R^{n_i}$, $u_{i,k}\in\R^{m_i}$, $w_{i,k}\in\R^{p_i}$, and $\sum_{i\in\cl{N}} n_i = n$, $\sum_{i\in\cl{N}} m_i = m$, $\sum_{i\in\cl{N}} p_i = p$. 
The following assumption is important in order to be able to partition the system parameters.

\myblue{
	\begin{assumption}\label{ass:decompose}
		The control input and the additive disturbance of the subsystems are decoupled, e.g., $u_{i,k}$ and $w_{i,k}$ only affect subsystem $i\in\cl{N}$ for all $k\in\cl{T}$. 
		Consider the state  and control input constraint sets $\mathcal{X}$ and $\mathcal{U}$ of large-scale system dynamics \eqref{dyn_LS} as defined in \eqref{feasible_set_LS} to be $\cl{X} = \prod_{i\in\cl{N}} \cl{X}_i$, and $\cl{U} = \prod_{i\in\cl{N}} \cl{U}_i$ such that  $\mathcal{X}_i$ and $\mathcal{U}_i$ for all subsystems $i\in\cl{N}$ can be in the following form:
		\begin{align*}
		\mathcal{X}_i := \{ {x} \in \mathbb{R}^{m_i} : G_i  {x} \leq g_i \}	, 
		\mathcal{U}_i := \{ {u} \in \mathbb{R}^{p_i} : H_i {u} \leq h_i  \},	
		\end{align*}
		where $G = \text{diag}_{i\in\mathcal{N}} (G_i)$, $H = \text{diag}_{i\in\mathcal{N}} (H_i)$, $g = \text{col}_{i\in\cl{N}}(g_i)$ and $h = \text{col}_{i\in\cl{N}}(h_i)$.
		Consider also the objective function of each subsystem $i\in\cl{N}$ in the following form:
		\begin{equation*}
		\begin{aligned}
		J_i(\bs{x}_i, \bs{u}_i) := \mathbb{E} \left[ \sum_{k=0}^{T-1} \left({x}_{i,k}^\top Q_i {x}_{i,k}+ {u}_{k}^\top R_i {u}_{i,k}\right) + {x}_{i,T}^\top P_i {x}_{i,T} \right] ,
		\end{aligned}
		\end{equation*}
		where $Q_i \in \mathbb{R}^{n_i \times n_i}_{\succeq 0}$, $R_i \in \mathbb{R}^{m_i \times m_i}_{\succ 0}$ such that $Q = \text{diag}_{i\in\cl{N}}(Q_i)$, and $R = \text{diag}_{i\in\cl{N}}(R_i)$. 
		Note that $\bs{x}_i = \text{col}_{k\in\cl{T}}(x_{i,k})$ and $\bs{u}_i = \text{col}_{k\in\cl{T}}(u_{i,k})$ such that $\bs{x} = \text{col}_{i\in\cl{N}}(\bs{x}_i)$ and $\bs{u} = \text{col}_{i\in\cl{N}}(\bs{u}_i)$.
	\end{assumption}
	
	It is important to note that under Assumption~\ref{ass:decompose} there are no direct coupling constraints between each subsystem $i\in\mathcal{N}$ and its neighboring subsystems $j\in\mathcal{N}_i$. 
	Instead, the subsystem $i\in\mathcal{N}$ is dynamically coupled with all its neighboring subsystems $j\in\mathcal{N}_i$ as it is presented in \eqref{dyn_sub}.
}

We refer to the additive disturbance $w_{i,k}$ as a private (local) uncertainty source of each subsystem $i\in\cl{N}$, since it is assumed that it affects only the subsystem $i\in\cl{N}$.
\myblue{Motivated by an application to Smart Thermal Grids (STGs) in \cite{rostampour2017probabilistic} where a network of agents are interconnected via an uncertain common resource pool between neighbors, we consider
	the uncertain variable $\delta_{k}$ as a common uncertainty source between all subsystems $i\in\cl{N}$.
	Such a STG application will also be presented in \Cref{ATES_STGs} as a second case study.}
Observe the fact that every common uncertain phenomenon can be considered as a local uncertain variable, e.g., the outside weather condition for neighboring buildings.
Therefore, we also consider to have ${\delta}_k = \text{col}_{i\in\cl{N}}(\delta_{i,k})$ and refer to both random variables $w_{i,k}$ and ${\delta}_{i,k}$ as local uncertainty sources.

We are now able to decompose the large-scale system matrices  
$B({\delta}_k) = \text{diag}_{i\in\mathcal{N}} (B_{i}({\delta}_{i,k}))\in \R^{n\times m}$, 
$C({\delta}_k) = \text{diag}_{i\in\mathcal{N}} (C_{i}({\delta}_{i,k}))\in \R^{n\times p}$,
and consider $A({\delta}_k) \in \R^{n\times n}$ in the following form:
\begin{equation*}
\begin{aligned}
A({\delta}_k) &= \begin{bmatrix}
A_{11}({\delta}_{1,k}) & \cdots & A_{1N}({\delta}_{1,k}) \\ 
\vdots & \ddots & \vdots \\
A_{N1}({\delta}_{N,k}) &  \cdots & A_{NN}({\delta}_{N,k})
\end{bmatrix} \,, 
\end{aligned}
\end{equation*}
where $A_{ij}({\delta}_{i,k}) \in \R^{n_i\times n_j}$, $B_{i}({\delta}_{i,k}) \in \R^{n_i\times m_i}$, and $C_i({\delta}_{i,k}) \in \R^{n_i\times p_i}$.
Define the set of neighboring subsystems of subsystem $i$ as follows:
\begin{align}\label{dyni}
\cl{N}_i = \left\{j\in\cl{N}\backslash \{i\} \ \big\vert \ A_{ij}({\delta}_{i,k}) \neq \bs{0} \right\} \ , 
\end{align}
where $\bs{0}$ denotes a matrix of all zeros with proper dimension.
Note that if subsystems are decoupled, they remain decoupled regardless of  the uncertainties ${\delta}_{i,k}\,$ for all $i\in\mathcal{N}$.
Consider now a large-scale network that consists of $N$ interconnected subsystems, and each subsystem can be described by 
an uncertain discrete-time linear time-invariant system with additive disturbance of the form 
\begin{align}\label{dyn_sub}
\begin{cases}
\, x_{i,k+1} &=  A_{ii}({\delta}_{i,k}) x_{i,k} + B_{i}({\delta}_{i,k}) u_{i,k} + q_{i,k}  \\
\, \ q_{i,k} &= \sum_{j\in\cl{N}_i} A_{ij}({\delta}_{i,k}) x_{j,k} + C_i({\delta}_{i,k}) w_{i,k} 
\end{cases} \ .
\end{align}
Following Assumption~\ref{ass:decompose}, one can consider a linear feedback gain matrix $K_i$ for each subsystem $i\in\cl{N}$ such that $K = \text{diag}_{i\in\cl{N}}(K_i)$.
Using $K_i$ in each subsystem, we assume that there exists $P_i$ for each subsystem $i\in\cl{N}$ such that $P = \text{diag}_{i\in\cl{N}}(P_i)$ preserves the condition in \eqref{feedback_law}.

\myblue{
	Define $\bs{\sigma}_i := \text{col}_{s\in\{1,\cdots,S_i\}}(\sigma_{i}^s)$ such that $\sigma_{i}^s:= \text{col}_{j\in\cl{N}_i}(\bs{z}_{ij}^s)$, where $\bs{z}_{ij}^s = \text{col}_{k\in\cl{T}}(z_{ij,k}^s)\in\R^{Tn_j}$, to be auxiliary decision variables, and using $\bs{v}_i = \text{col}_{k\in\cl{T}}(v_{i,k})\in\R^{Tm_i}$ such that $\bs{v} = \text{col}_{i\in\cl{N}}(\bs{v}_i)$, we are now able to reformulate the optimization problem in \eqref{scp_ls} into the following decomposable finite horizon scenario program: 
	\begin{subequations}
		\label{scp_sub}
		\begin{align}
		\underset{
			\begin{subarray}{c}
			\{\bs{\sigma}_i,\bs{v}_i\}_{i\in\mathcal{N}}
			\end{subarray}
		}{\min} & \sum_{i\in\mathcal{N}}^{}J_i(\bs{x}_i, \bs{u}_i)  \\
		\text{s.t.} \ \quad &  {x}_{i,k+1}^{s} = A_{ii}(\delta_{i,k}^{s}) {x}_{i,k}^{s} + B_{i}(\delta_{i,k}^{s}) {u}_{i,k}^{s} + {q}_{i,k}^{s} \,, \label{scp_sub_dyn} \\
		& q_{i,k}^{s} = \sum\nolimits_{j\in\cl{N}_i} A_{ij}({\delta}_{i,k}^{s}) z_{ij,k}^{s} + C_i({\delta}_{i,k}^{s}) w_{i,k}^{s}, \label{scp_sub_qi} \\
		& z_{ij,k}^{s} = x_{j,k}^{s} \,,\, \forall\bs{x}_j^{s}\in\cl{S}_{\bs{x}_j}\,, \forall j\in\cl{N}_i\,, \label{cons_consensus} \\
		&  {x}_{i,k+\ell}^{s} \in \cl{X}_i, \ell \in \cl{T}_+, \forall\bs{w}^{s}_i \in {\cl{S}}_{\bs{w}_i}, \forall\bs{\delta}_i^{s} \in {\cl{S}}_{\bs{\delta}_i}, \label{scp_cons_sub} \\
		&  {u}_{i,k}^{s} = K_{i} {x}_{i,k}^{s} + {v}_{i,k} \in \cl{U}_i \,,\, \forall k\in\cl{T} \,, \forall i\in\mathcal{N}\, \label{cons_u_sub}
		\end{align}
	\end{subequations}
	where $\bs{w}_i = \text{col}_{k\in\cl{T}}(w_{i,k})\in\cl{W}_i$ and $\bs{\delta}_i = \text{col}_{k\in\cl{T}}(\delta_{i,k})\in\Delta_i$
	such that $\cl{W} =  \prod_{i\in\cl{N}}\cl{W}_{i}$ and $\Delta =  \prod_{i\in\cl{N}}\Delta_{i}$. 
	The sets $\cl{S}_{\bs{w}_i}:= \{\bs{w}^{1}_i, \cdots, \bs{w}^{S_i}_i\}\in\cl{W}^{S_i}_i$ and $\cl{S}_{\bs{\delta}_i}:= \{\bs{\delta}^{1}_i, \cdots, \bs{\delta}^{S_i}_i\}\in\Delta^{S_i}_i$ denote sets of given finite samples (scenarios) of disturbance and uncertainties in each subsystem $i\in\cl{N}$, such that $\cl{S}_{\bs{w}} =  \prod_{i\in\cl{N}}\cl{S}_{\bs{w}_i}$ and $\cl{S}_{\bs{\delta}} =  \prod_{i\in\cl{N}}\cl{S}_{\bs{\delta}_i}$.
	It is important to highlight that the proposed optimization problem in \eqref{scp_sub} is decomposable into $N$ subproblems with coupling only through the equality constraints \eqref{cons_consensus} between neighboring agents.
	The auxiliary variables $z_{ij,k}^{s}$ are defined between  agent $i\in\cl{N}$ and its neighboring agent $j\in\cl{N}_j$ to represent the local estimation of each scenario information of the neighbor's state variable $\bs{x}_j^{s} = \text{col}_{k\in\cl{T}}(x_{j,k}^{s}) \in \mathcal{X}_j^{T}:= \mathbb{X}_j$ at each sampling time  $k\in\cl{T}$.
	By taking into consideration that the interaction dynamics model $A_{ij}({\delta}_{i,k}^{s})$ by each neighboring agent $j\in\mathcal{N}_i$ is available for agent $i\in\mathcal{N}$, the only information that subsystem $i\in\cl{N}$ needs is a set of scenarios of the state variable 
	$\cl{S}_{\bs{x}_j} :=\{\bs{x}_j^{1}, \cdots, \bs{x}_j^{S_i}\}\in \mathbb{X}_j^{S_i}$ at each
	sampling time $k\in\cl{T}$ from all its neighboring agents $j\in\cl{N}_i$.
	Note that the cardinality of $\cl{S}_{\bs{x}_j}$ is $S_i$, which is assigned by the agent $i\in\mathcal{N}$ for all neighboring agents $j\in\cl{N}_i$, and should be sent by the neighboring agents at each
	sampling time $k\in\cl{T}$, and this is called as the \textit{hard communication scheme} between the neighboring agents.
}

In the following proposition, we provide a connection between the proposed optimization problem in \eqref{scp_sub} and the optimization problem in \eqref{scp_ls}.

\myblue{
	\begin{proposition}\label{propos:decomp_exact}
		Given Assumption~\ref{ass:decompose} and the block-diagonal structure for the state-feedback controller $K = \text{diag}_{i\in\cl{N}}(K_i)$ for the large-scale system dynamics \eqref{dyn_LS}, the optimization problem in \eqref{scp_sub} is an exact decomposition of the optimization problem in \eqref{scp_ls}, such that the optimal objective value of the decomposed problem \eqref{scp_sub} is equal to
		the optimal objective of the centralized problem \eqref{scp_ls}.
	\end{proposition}
}

The proof is provided in the Appendix. \QEDB

\myblue{
	An important key feature of the proposed decomposable scenario program in \eqref{scp_sub} compared to the centralized scenario problem in \eqref{scp_ls} is as follows.
	Such a decomposition yields a reduction of computation time complexity of the centralized scenario program compared to the decomposable scenario program \eqref{scp_sub}.
	This however requires more communication between each subsystem, since at each
	sampling time $k\in\cl{T}$ each agent $i\in\cl{N}$ should send to and receive from all the neighboring agents $j\in\cl{N}_i$ a set of their state variable scenarios $\cl{S}_{\bs{x}_i}$ and $\cl{S}_{\bs{x}_j}$, respectively\ifshort.\else, as it is shown in Fig.~\ref{fig:dsmpc}.\fi
}

\begin{remark}\label{ccp_sub}
	The proposed constraint \eqref{scp_cons_sub} represents an approximation of the following chance constraint on the state of each subsystem $i\in\cl{N}$:
	\begin{align}\label{ccp_cons_sub}
	\Pb [\, {x}_{i,k+\ell} \in \cl{X}_i \,,\, \ell \in \cl{T}_+\,] \geq 1 - \varepsilon_i \ ,
	\end{align}
	where $\varepsilon_i\in(0,1)$ is the admissible state constraint violation parameter of each subsystem \eqref{dyn_sub}.
	One can also consider $\alpha_i = 1-\varepsilon_i$ as the desired level of state feasibility parameter of each subsystem \eqref{dyn_sub}.
\end{remark}

The following theorem can be considered as the main result of this section to quantify the robustness of the solutions obtained by \eqref{scp_sub}.

\begin{theorem}\label{scenario2_thm}
	Let $\varepsilon_i, \beta_i \in (0,1)$ be chosen such that $\varepsilon = \sum_{i\in\cl{N}} \varepsilon_i \in (0,1)$, $\beta = \sum_{i\in\cl{N}} \beta_i \in (0,1)$, and $N_{s_i} \geq \mathsf{N}(\varepsilon_i, \beta_i, Tm_i)$  for all subsystem $i\in\cl{N}$.
	If $\bs{v}^* = \text{col}_{i\in\cl{N}}(\bs{v}_i^*)$, the collection of the optimizers of problem \eqref{scp_sub} for all subsystem $i\in\cl{N}$, is applied to the discrete-time dynamical system \eqref{dyn_LS} for a finite horizon of length $T$, then, with at least confidence $1-\beta$, the original constraint \eqref{ccp_ls} is satisfied for all $k\in\cl{T}$ with probability more than $1-\varepsilon$.
\end{theorem}

The proof is provided in the Appendix. \QEDB

The interpretation of \Cref{scenario2_thm} is as follows.
In the proposed decomposable scenario program \eqref{scp_sub}, each subsystem $i\in\cl{N}$ can have a desired level of constraint violation $\varepsilon_i$ and a desired level of confidence level $1-\beta_i$.
To keep the robustness level of the collection of solutions in a probabilistic sense \eqref{ccp_ls} for the discrete-time dynamical system \eqref{dyn_LS}, these choices have to follow a certain design rule, e.g. $\varepsilon = \sum_{i\in\cl{N}} \varepsilon_i \in (0,1)$  and $\beta = \sum_{i\in\cl{N}} \beta_i \in (0,1)$. 
This yields a fixed $\varepsilon\,,\, \beta$ for the large-scale system \eqref{dyn_LS} and the individual $\varepsilon_i\,,\, \beta_i$ for each subsystem $i\in\cl{M}$. 
It is important to mention that in order to maintain the violation level for the large-scale system with many partitions, i.e. $|\mathcal{N}| = N \uparrow\,,$ the violation level of individual agents needs to decrease, i.e. $\varepsilon_i \downarrow\,,$ which may lead to conservative results for each subsystem, since the number of required samples needs to increase, i.e. $S_i\uparrow$.
Addressing such a limitation is subject of ongoing research work.

\myblue{
	The following corollary is a direct result of \Cref{scenario2_thm}.
	\begin{corollary} 
		If the optimal solution $\bs{v}_i^*$ for each agent $i\in\mathcal{N}$ obtained via the proposed decomposable problem \eqref{scp_sub} is applied to the agent dynamical system \eqref{dyn_sub} for a finite horizon of length $T$, then, with at least confidence $1-\beta_i$, the original constraint \eqref{ccp_cons_sub} is satisfied for all $k\in\cl{T}$ with probability more than $1-\varepsilon_i$.
	\end{corollary}
}

\ifshort
\else

\begin{figure}
	\centering
	\includegraphics[width=0.65\textwidth]{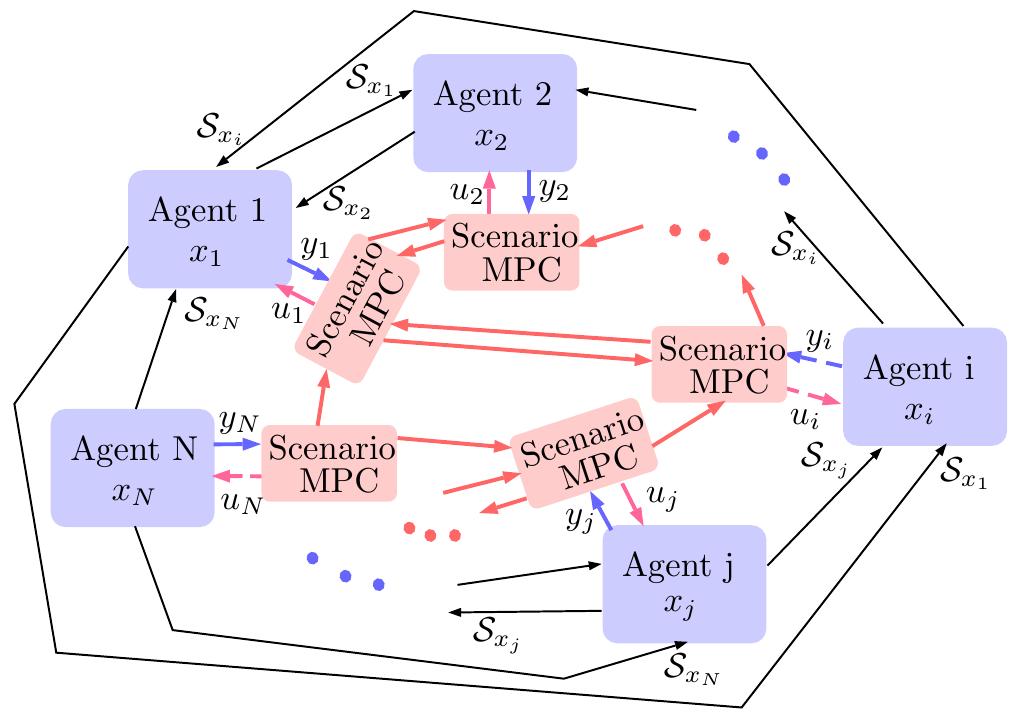}
	\caption{
		\myblue{
			Distributed scenario MPC which corresponds to the problem \eqref{scp_sub}.
			The measurement variable $y_i$ is the full state information $x_i$ which is sent to the distributed controllers (local scenario MPC) for all agent $i\in\mathcal{N}$. The red lines present the communication needed between the neighboring controllers following the proposed distributed scenario exchange scheme (\Cref{alg:distributed scenario exchange ADMM}).}
	}
	\label{fig:dsmpc}
\end{figure}

To summarize this section, we present a network of agents that are dynamically coupled and with their own local scenario MPC is depicted in Fig.~\ref{fig:dsmpc}. 
\myblue{
	In the next section, we will develop a distributed framework based on the ADMM algorithm to solve the proposed optimization problem \eqref{scp_sub} using the so-called hard communication scheme.
}

\fi

\myblue{
\section{Distributed Scenario MPC}\label{dist_smpc}

In this section, we continue by developing a distributed framework using ADMM to solve the proposed decomposable formulation in \eqref{scp_sub}. 
It has been proven that ADMM for this type of problem converges linearly \cite{shi2014linear}.
We follow a similar approach as in \cite{rostampour2017distributedTPS,boyd2011distributed} to solve the decomposable problem \eqref{scp_sub} in a distributed manner by extending to the scenario program.

Define the local feasibility set denoted with $\mathcal{L}_i$ for all agents $i \in \mathcal{N}$:
\begin{align*} 
\mathcal{L}_i := \Big\{ &\{ \bs{\sigma}_i \,, \bs{v}_i \} \, \Big| \, {u}_{i,k}^{s} = K_{i} {x}_{i,k}^{s} + {v}_{i,k} \in \cl{U}_i, \forall k\in\cl{T}, \\ 
&{x}_{i,k+1}^{s} = A_{ii}(\delta_{i,k}^{s}) {x}_{i,k}^{s} + B_{i}(\delta_{i,k}^{s}) {u}_{i,k}^{s} + {q}_{i,k}^{s},  \\
& q_{i,k}^{s} = \sum\nolimits_{j\in\cl{N}_i} A_{ij}({\delta}_{i,k}^{s}) z_{ij,k}^{s} + C_i({\delta}_{i,k}^{s}) w_{i,k}^{s}, \\ 
& {x}_{i,k+\ell}^{s} \in \cl{X}_i, \ell \in \cl{T}_+, \forall\bs{w}^{s}_i \in {\cl{S}}_{\bs{w}_i}, \forall\bs{\delta}_i^{s} \in {\cl{S}}_{\bs{\delta}_i}\ \ \Big\}\,, 
\end{align*}
where all the system parameters together with the local scenario sets, ${\cl{S}}_{\bs{w}_i}, {\cl{S}}_{\bs{\delta}_i}$, are fixed and available locally for each agent $i\in\mathcal{N}$ at each sampling time $k\in\mathcal{T}$. 
We are now able to define the augmented Lagrangian function as follows:
\begin{align*}
L(\{\bs{\sigma}_i \,, \bs{v}_i \}_{i\in\mathcal{N}}) &:= \sum_{i \in \mathcal{N}}  \Bigg( J_i(\bs{x}_i, \bs{u}_i)
+ I_{\mathcal{L}_i}(\{ \bs{\sigma}_i, \bs{v}_i \}) \ + 
\sum_{j \in \mathcal{N}_i}  \sum_{s = 1}^{S_i} \Big( \frac{\mu}{2} \Big\| \bs{z}^s_{ij} - \bs{x}^s_{j} + \frac{\Lambda_{ij}^s}{\mu} \Big\|^2_2 + \frac{1}{2\mu} \Big\| \Lambda_{ij}^s \Big\|_{2}^{2} \Big) \Bigg),
\end{align*}
where $I_{\mathcal{L}_i}(\{ \bs{\sigma}_i, \bs{v}_i \})$ is a convex indicator function for the local constraint set $\mathcal{L}_i$ that maps to infinity if one of the constraints is violated, and to zero otherwise.
Step size $\mu$ is a fixed constant, and multipliers $\Lambda_{ij}^s \in \R^{n_j}$ are introduced for every coupling constraint of each scenario. 
The indicator function makes sure all local constraints are satisfied, and the squared norm penalty term forces constraint~\eqref{cons_consensus} to be satisfied.

We now describe the steps of the ADMM algorithm as follows:

\subsubsection*{(1) Update primal variables} 
The multipliers $\Lambda_{ij}^s$ and the local estimation of neighboring state scenario $\bs{z}_{ij}^s$ for all scenarios $s=1,\cdots, S_i$ and for all the neighbors $j\in\mathcal{N}_i$ are fixed at their value of the previous iteration.
Consider the following minimization problem for all $i \in \mathcal{N}$:
\begin{align}\label{eq:admm multi-area primal update}
\bs{v}_i^{(l+1)} = \arg&\min_{\bs{v}_i} \Bigg\{ J_i(\bs{x}_i, \bs{u}_i)
+ I_{\mathcal{L}_i}(\{ \bs{\sigma}_i^{(l)}, \bs{v}_i \}) \ +  \sum_{j \in \mathcal{N}_i}  \sum_{s = 1}^{S_j} \Big( \frac{\mu}{2} \Big\| \bs{z}^{s,(l)}_{ji} - {\bs{x}}^{s}_{i} + \frac{\Lambda_{ji}^{s,(l)}}{\mu} \Big\|^2_2 \Big)   \Bigg\} \,.
\end{align}
It is important to note that the augmented part of the aforementioned step in each agent $i\in\mathcal{N}$ refers to the coupling constraints from the neighboring agent $j\in\mathcal{N}_i$ point of view and $\bs{z}^{s,(l)}_{ji}$ is the estimation of neighboring agent $j\in\mathcal{N}_i$ from the state variable scenario of agent $i\in\mathcal{N}$ for all $s \in\{1, \cdots, S_j\}$.
Since the minimization is only in $\bs{v}_i$, all terms of $\Lambda_{ij}^s$ drop out.
This results in $|\mathcal{N}| = N$ separate convex programs.

\subsubsection*{(2) Update and exchange the set of state scenarios} 
The resulting primarily decision variables $\bs{v}_i^{(l+1)}$ for all $i \in \mathcal{N}$ are used to first update and then exchange the set of state scenarios.
Note that each agent only needs to communicate the updated local state scenarios set to all its neighboring agents and receive the updated set of all its neighboring state scenarios.
Consider now the following steps which we refer to as the Scenario Updating Steps (SUS) for a generic agent $i\in\mathcal{N}$ as follows:
\begin{align}\label{SUS}
\begin{cases}
{x}_{i,k+1}^{s,(l+1)} = A_{ii}(\delta_{i,k}^{s}) {x}_{i,k}^{s,(l+1)} + B_{i}(\delta_{i,k}^{s}) {u}_{i,k}^{s,(l+1)} + {q}_{i,k}^{s,(l)}  \\ 
q_{i,k}^{s,(l)} = \sum\nolimits_{j\in\cl{N}_i} A_{ij}({\delta}_{i,k}^{s}) z_{ij,k}^{s,(l)} + C_i({\delta}_{i,k}^{s}) w_{i,k}^{s} \\
{u}_{i,k}^{s,(l+1)} = K_{i} {x}_{i,k}^{s,(l+1)} + {v}_{i,k}^{(l+1)}, \forall k\in\cl{T} 
\end{cases}\,.
\end{align}
By repeating SUS for all local (private) scenarios, $\forall\bs{w}^{s}_i \in {\cl{S}}_{\bs{w}_i}$ and $\forall\bs{\delta}_i^{s} \in {\cl{S}}_{\bs{\delta}_i}$, 
the agent $i\in\mathcal{N}$ should first update the set of its local state scenarios $\bs{x}^{s,(l+1)}_{i}\in\mathcal{S}^{(l+1)}_{\bs{x}_i}$ and send it to all its neighboring agents.
Then, the agent $i\in\mathcal{N}$ should receive the updated set of all its neighboring state scenarios $\mathcal{S}^{(l+1)}_{\bs{x}_j}$.
It is important to note that the set of local scenarios, ${\cl{S}}_{\bs{w}_i}$ and ${\cl{S}}_{\bs{\delta}_i}$, is fixed for all the SUS steps and iteration steps, i.e., $l\in\N_+$.

\subsubsection*{(3) Update and exchange the  state estimation variables} 

Using the updated set of all neighboring state scenarios  $\mathcal{S}^{(l+1)}_{\bs{x}_j}$ and the updated primal variables, consider the following projection problem for all $i \in \mathcal{N}$:
\begin{align}\label{eq:admm multi-area shared update}
\bs{\sigma}_i^{(l+1)} = \arg&\min_{\bs{\sigma}_i} \Bigg\{I_{\mathcal{L}_i}(\{ \bs{\sigma}_i^{}, \bs{v}_i^{(l+1)} \}) \ + 
\sum_{j \in \mathcal{N}_i}  \sum_{s = 1}^{S_i} \Big( \frac{\mu}{2} \Big\| \bs{z}^{s}_{ij} - {\bs{x}}^{s,(l+1)}_{j} + \frac{\Lambda_{ij}^{s,(l)}}{\mu} \Big\|^2_2 \Big)   \Bigg\} \,,
\end{align}
%
%
Since the projection is only in $\bs{\sigma}_i$, all terms of $\Lambda_{ij}^s$ drop out and this results in $|\mathcal{N}| = N$ separate convex programs.
Each agent $i \in \mathcal{N}$ should now send the updated local estimation of the neighboring state variable $\bs{z}^{s,(l+1)}_{ij}$ to its neighbor $j \in \mathcal{N}_i$.

\subsubsection*{(3) Update and exchange the  multiplier variables} 
The multipliers are updated as follows $\forall i \in \mathcal{N}, \forall j \in \mathcal{N}_i$, and $s = 1, \cdots, S_i$:
\begin{align}\label{eq:admm multi-area dual update}
\Lambda_{ij}^{s,(l+1)} &= \Lambda_{ij}^{s,(l)} + \mu \left( \bs{z}^{s,(l+1)}_{ij} - \bs{x}^{s,(l+1)}_{j} \right) \,.
\end{align}
Notice that no information exchange is needed for the update of the multiplier, since agents receive the updated set of all their neighboring state scenarios $\bs{x}^{s,(l+1)}_{j}\in\mathcal{S}^{(l+1)}_{\bs{x}_j}$.
As the final step, each agent $i \in \mathcal{N}$ should receive the updated multiplier variable $\Lambda_{ji}^{s,(l+1)}$ from all its neighboring agents $j \in \mathcal{N}_i$ for all the scenarios $s = 1, \cdots, S_j$.

In \Cref{alg:distributed scenario exchange ADMM}, we summarize the proposed distributed scenario exchange framework using ADMM algorithm to illustrate the calculation and communication steps for each agent $i\in\mathcal{N}$. 
Each agent needs to solve two small-scale scenario programs in \Cref{line:maadmm primal update} and \Cref{line:maadmm shared update} at each iteration that can be considered as the highest computational cost in the proposed algorithm.
In \Cref{line:maadmm simulation} each agent first updates $\mathcal{S}_{\bs{x}_i}^{(l+1)}$ using SUS, that works with only simple operations, e.g., addition, subtraction and scaling, and then, it broadcasts and receives the updated local $\mathcal{S}_{\bs{x}_i}^{(l+1)}$ and $\mathcal{S}_{\bs{x}_j}^{(l+1)}$ to and from the neighbors, respectively.
Finally, each agent in \Cref{line:maadmm broadcast shared} and \Cref{line:maadmm broadcast multipliers} broadcasts the updated local estimation from the neighboring agents state variable scenarios and the related multipliers.

\begin{algorithm}[t]
	\caption[Distributed Scenario Exchange Algorithm]{Distributed Scenario Exchange Algorithm}
	\label{alg:distributed scenario exchange ADMM}
	\begin{algorithmic}[1]
		\ForAll {$i \in \mathcal{N}$}
		\State \textbf{initialize}: ${x}_{i,0} \in \cl{X}_i$, $l = 0$, $\Lambda_{ij}^{s,(l)}$, $\bs{z}_{ij}^{s,(l)}$, and $\bs{z}_{ji}^{s,(l)}$
		\Statex \ \quad $\forall s\in\{1,\cdots, S_i\}$, and $\forall j \in \mathcal{N}_i$ 
		\While {$\eta^{(l)}_i > \epsilon^{\text{des}}_i$ using \eqref{convergence criteria ADMM}}
		\State \label{line:maadmm primal update}
		\textbf{Update} $\bs{v}_i^{(l+1)}$ using \eqref{eq:admm multi-area primal update} 	
		\State \textbf{Update} $\mathcal{S}_{\bs{x}_i}^{(l+1)}$ using SUS procedure \eqref{SUS} \label{line:maadmm simulation} 
		\State \label{line:maadmm broadcast} 
		\textbf{Broadcast} $\mathcal{S}_{\bs{x}_i}^{(l+1)}$ to all $j \in \mathcal{N}_i$
		\State \label{line:maadmm receive} 
		\textbf{Receive} $\mathcal{S}_{\bs{x}_j}^{(l+1)}$ from all $j \in \mathcal{N}_i$
		\State \label{line:maadmm shared update}
		\textbf{Update} $\bs{z}_{ij}^{s,(l+1)}$ using \eqref{eq:admm multi-area shared update} for all $j \in \mathcal{N}_i$ and 
		\Statex \qquad \quad for all $ s\in\{1,\cdots, S_i\}$, 
		\State \label{line:maadmm broadcast shared} 
		\textbf{Broadcast} $\bs{z}_{ij}^{s,(l+1)}$ to all $j \in \mathcal{N}_i$
		\State \textbf{Update} $\Lambda_{ij}^{s,(l+1)}$ using \eqref{eq:admm multi-area dual update} for all $j \in \mathcal{N}_i$ and  
		\Statex \qquad\quad for all $s\in\{1,\cdots, S_i\}$
		\State \label{line:maadmm broadcast multipliers} 
		\textbf{Broadcast} $\Lambda_{ij}^{s,(l+1)}$ to all $j \in \mathcal{N}_i$ and  
		\Statex \qquad\quad for all $s\in\{1,\cdots, S_i\}$
		\State \textbf{set} \ $l \gets l+1$
		\EndWhile
		\State \textbf{output:} $\bs{v}^*_i = \bs{v}^{l+1}_i$ 
		\EndFor 
	\end{algorithmic}
\end{algorithm}

We define the agreement on the local estimation from the state variable scenario of the neighbors $\bs{z}_{ij}^s$ as the convergence criteria of the proposed ADMM algorithm for each agent $i \in \mathcal{N}$ as follows:
\begin{align}\label{convergence criteria ADMM}
\eta^{(l)}_i =  \sum_{j \in \mathcal{N}_i} \sum_{s = 1}^{S_i} \ \left( \frac{\mu}{2} \Big\| \bs{z}^{s,(l)}_{ij} - \bs{x}^{s,(l)}_{j} \Big\|_{2}^{2}  \right) \ , 
\end{align}
and the overall convergence criteria is $\eta^{(l)}:=\sum_{i \in \mathcal{N}} \eta^{(l)}_i$.
If the residual sequence $\left\{\eta^{(l)}\right\}_{l=1}^{+\infty}$ is sufficiently small, less than or equal to $\sum_{i \in \mathcal{N}}\epsilon^{\text{des}}_i$ where $\epsilon^{\text{des}}_i$ is the desired level of convergence error for each agent $i\in \mathcal{N}$, then all neighboring agents have reached a consistent set of state variable scenarios. 
In light of Assumption~\ref{ass_scp_ls_feasible}, we can now provide convergence of \Cref{alg:distributed scenario exchange ADMM} based on the results in  \cite{he2015non}.
\begin{theorem}\label{thm_admm_opf}
	Assume that Slater's condition  \cite{slater1959lagrange} holds for the decomposable scenario optimization problem \eqref{scp_sub}, and consider the iterative steps given in \Cref{alg:distributed scenario exchange ADMM}. 
	Then the following statements hold:
	\begin{itemize}
		\item The residual sequence $\{\eta^{(k)}\}_{k=0}^{+\infty}$ tends to $0$ in a non-increasing way as $l$ goes to $+\infty$, and consequently, for all $j\in\mathcal{N}_j\,,$ and each $i\in\mathcal{N}$: 
		\begin{align*}
		\bs{z}^{s,(+\infty)}_{ij} = \bs{x}^{s,(+\infty)}_{j}\,.
		\end{align*}
		\item The sequence $\{\bs{v}_i^{(l)}\}_{\forall i\in\mathcal{N}}$ generated by \Cref{alg:distributed scenario exchange ADMM} converges to an optimal solution $\{\bs{v}_i^{*}\}_{\forall i\in\mathcal{N}}$ of the decomposable scenario program \eqref{scp_sub} as $l$ tends to $+\infty$.	
	\end{itemize}
\end{theorem}
\textit{Proof.}
The theorem follows from \cite{he2015non} that studies the  convergence  of a standard ADMM  problem.  
The  details are omitted for brevity. 
\QEDB

\begin{remark}
	The proposed distributed scenario exchange algorithm uses a Gauss-Siedel update on the primal variables and the set of state variable scenarios, after which the multiplier variables are updated.
	Since either the primal or the set of state variable scenarios are fixed in the Gauss-Siedel steps, the problem can be distributed for all the agents.
	The main advantage to distribute a large-scale scenario optimization problem \eqref{scp_ls} is the ability of finding local solutions for each agent based on the information received in the previous iteration. 
	Such calculations can therefore be carried out in parallel.
	Although an actual parallel implementation is outside the scope of this work, it is important to mention that the proposed algorithm is amenable to such an implementation.
	ADMM algorithms typically need a large number of iterations to converge to high accuracy, so the local agent problems need to be solved many times before finding a good enough solution.
	Thus, the ADMM approach without parallelization might not be the quickest method to solve the large-scale scenario program.
	Such a distributed framework is advantageous especially when the global scenario optimization problem is too large and cannot be solved within polynomial time due to the curse of dimensionality or memory limitations and computational constraints.
\end{remark}

We next present the overall steps needed to execute the proposed distributed scenario MPC.
\Cref{DSMPC} summarizes all the steps such that after decomposing the large-scale dynamical system \eqref{dyn_LS} and determining the index sets of the neighboring agents, each agent first starts to generate some scenarios of its local uncertainty sources in order to approximate its cost function and achieve the feasibility of its chance constraint with high confidence level following \Cref{scenario2_thm}. 
Next, all the agents solve their own local problem using the proposed distributed scenario exchange scheme (\Cref{alg:distributed scenario exchange ADMM}) to obtain their local optimal solution. 
To execute \Cref{alg:distributed scenario exchange ADMM}, it is assumed that the feedback control gain matrices $K_i$ for all agents $i\in\cl{N}$ are given \eqref{feedback_law}, and the coupling terms $A_{ij}(\delta_{i,k})$ are known between each agent $i\in\cl{N}$ and its neighboring agents $j\in\cl{N}_i$.
Using the obtained solution via  \Cref{alg:distributed scenario exchange ADMM}, each agent generates a new set of local state variable scenarios to send and receive to and from all its neighboring agents, respectively. 
Finally, the first optimal control input is applied to the real system and the new state variables are measured, and proceed in the classical receding horizon fashion.

\begin{algorithm}[t]
	\caption{Distributed Scenario MPC (DSMPC)}\label{DSMPC}
	\small
	\begin{algorithmic}[1]
		\State \textbf{Decompose} the large-scale dynamical system \eqref{dyn_LS} into $N$ agents as the proposed form in \eqref{dyn_sub}
		\State \textbf{Determine} the index set of neighboring agents $\cl{N_i}$ for each agent $i\in\cl{N}$
		\ForAll{agent $i\in\cl{N}$}
		\State \textbf{fix} initial state ${x}_{i,0} \in \cl{X}_i$, $\varepsilon_i\in(0,1)$, and  $\beta_i\in(0,1)$ 
		\Statex \quad\,\ such that $\varepsilon = \sum\nolimits_{i\in\cl{N}}\varepsilon_i\in(0,1) \,,\, \beta = \sum\nolimits_{i\in\cl{N}}\beta_i\in(0,1)$ \label{start dsmpc}
		\State \textbf{determine} $\bar{S}_i\in(0,+\infty)$ to approximate the objective 
		\Statex \quad\,\ function, and $S_{i}$ following \Cref{scenario2_thm} to approximate 
		\Statex \quad\,\ the chance constraint \eqref{ccp_cons_sub} in an equivalent sense
		\State \textbf{generate} $\bar{S}_i$, $S_{i}$ scenarios of $\bs{w}_i$, $\bs{\delta}_i$ to determine the sets 
		\Statex \quad\,\ of $\bar{\cl{S}}_{\bs{w}_i,\bs{\delta}_i}$ and $\cl{S}_{\bs{w}_i}$, $\cl{S}_{\bs{\delta}_i}$ 
		\State {\textbf{solve} the proposed optimization problem in \eqref{scp_sub} using
			\Statex  \quad\,\   the proposed distributed scenario exchange algorithm
			\Statex  \quad\,\  (\Cref{alg:distributed scenario exchange ADMM}) and determine  $\bs{v}^*_i$}
		\State \textbf{generate} ${S}_{i}$ scenarios of $\bs{x}_i$ using SUS procedure \eqref{SUS} 
		\Statex  \quad\,\  and taking into consideration $\bs{v}^*_i$ and $\cl{S}_{\bs{w}_i}$, $\cl{S}_{\bs{\delta}_i}$
		\State \textbf{send} the set ${\cl{S}}_{\bs{x}_i}$ to all neighboring agents $j\in\cl{N_i}$
		\State \textbf{receive} the sets ${\cl{S}}_{\bs{x}_j}$ from all neighboring agents $j\in\cl{N_i}$
		\State \textbf{apply} the first input of solution $u_{i,k}^* = K_{i} x_{i,k} + v^*_{i,k}$ into 
		\Statex \quad\,\ the uncertain subsystem \eqref{dyn_sub} 
		\State \textbf{measure} the state and substitute it as the initial state of 
		\Statex \quad\,\ the next step ${x}_{i,0}$ 
		\State \textbf{set} \ $k \gets k+1$ 
		\State \textbf{go to} Step \eqref{start dsmpc}
		\EndFor
	\end{algorithmic}
\end{algorithm}

\begin{remark}
	The proposed decomposable scenario program in \eqref{scp_sub} is a general deterministic optimization problem for each agent that are coupled via \eqref{cons_consensus}. Therefore, the proposed technique to solve \eqref{scp_sub} in \Cref{alg:distributed scenario exchange ADMM}, namely the distributed scenario exchange scheme, can be easily considered as a general solution method to solve the agents' problem with the nonempty intersection between their local feasible sets. Such a case where agents are only coupled via coupling constraints has been presented in \Cref{ATES_STGs} as a second case study. 
\end{remark}

\ifshortt
\else

To conclude, in this section we developed a solution method to solve the proposed decomposable scenario program in \eqref{scp_sub} using ADMM algorithm. However, such a technique requires many iteration steps at each time step, (see Step 7 in \Cref{DSMPC} to execute \Cref{alg:distributed scenario exchange ADMM} at each sampling time $k\in\mathcal{T}$), in order for the neighboring agents to agree on their possible set of state variable trajectories.
To overcome such an expensive communication requirement, in the next section, we will propose a novel inter-agent information exchange scheme to provide a more flexible framework to exchange information between agents, which subsequently yields a PnP distributed scenario MPC framework.

\fi

}

\section{Information Exchange Scheme}\label{inter_agent}

\ifshortt
\else

In this section, we first describe the information exchange between agents and then propose a set-based information exchange scheme which will be referred to as a {\em soft communication protocol} later in this section.
We finally provide the theoretical results to quantify robustness of the proposed information exchange scheme between neighboring agents. 

\fi

\ifshortt

When the proposed distributed scenario MPC framework (\Cref{DSMPC}) is applied to the large-scale scenario program \eqref{scp_ls}, all neighboring agent $j\in\cl{N}_i$ of the agent $i\in\cl{N}$ should send a set of scenarios of the state variable $\cl{S}_{\bs{x}_j} :=\{\bs{x}_j^{1}, \cdots, \bs{x}_j^{S_i}\}\in \mathbb{X}_j^{S_i}$ to agent $i$ at each sampling time $k\in\cl{T}$ following the proposed distributed scenario exchange scheme in \Cref{alg:distributed scenario exchange ADMM}.
Based on ADMM in \Cref{alg:distributed scenario exchange ADMM}, it may require a large number of iterations to achieve an agreement between neighboring agents on the consistency of exchange scenarios, which may also turn out to be too costly in terms of required communication bandwidth.
To address this shortcoming, we propose a set-based information exchange scheme which is referred to as a {\em soft communication protocol}.

\else

When the proposed distributed scenario MPC framework (\Cref{DSMPC}) is applied to the large-scale scenario program \eqref{scp_ls}, all neighboring agent $j\in\cl{N}_i$ of the agent $i\in\cl{N}$ should send a set of scenarios of the state variable $\cl{S}_{\bs{x}_j} :=\{\bs{x}_j^{1}, \cdots, \bs{x}_j^{S_i}\}\in \mathbb{X}_j^{S_i}$ to agent $i$ at each sampling time $k\in\cl{T}$ following the proposed distributed scenario exchange scheme in \Cref{alg:distributed scenario exchange ADMM}.
It is of interest to address the issue of how an agent $j\in\cl{N}_i$ can send the contents of $\cl{S}_{\bs{x}_j}$ to the agent $i\in\cl{N}$.

We propose the following two schemes:
1) Following our proposed setup in \eqref{DSMPC} to achieve a probabilistic guarantee for the obtained solution, agent $i\in\cl{N}$ requests from its neighboring agents to send the complete set of data $\cl{S}_{\bs{x}_j}$, element by element such that the number of required samples $S_i$, is chosen according to \Cref{scenario2_thm} in order to have a given probabilistic guarantee for the optimizer $\bs{v}_i^*$. 
We refer to this scheme as a {\em hard communication protocol} between agents following the proposed distributed scenario exchange scheme in \Cref{alg:distributed scenario exchange ADMM}.
Its advantage is that it is simple and transmits exactly the contents of $\cl{S}_{\bs{x}_j}$, but due to possibly high values of $S_i$, it may turn out to be too costly in terms of required communication bandwidth.
2) To address this shortcoming, we propose another scheme, where agent $j\in\cl{N}_i$ sends instead a suitable parametrization of a set that contains all the possible values of data with a desired level of probability ({\em the level of reliability}) $\tilde{\alpha}_j$. 
By considering a simple family of sets, for instance boxes in $\R^{n_j}$, communication cost can be kept down at reasonable levels.
We refer to this scheme as a {\em soft communication protocol} between agents\ifshort.\else (see Fig.~\ref{fig:inter_agent}).\fi
Such a scheme may be understood as a cascaded scenario scheme similar to the one in \cite{margellos2015connection}, where a sufficient number of scenarios was determined in order to establish probabilistic feasibility for two cascaded scenario programs subject to a similar source of uncertainty. 
Our soft communication setting however differs from \cite{margellos2015connection}, since each agent is subject to its own uncertainty source.
We aim to incorporate the reliability notion of such a soft communication scheme into the feasibility bound of each agent, in addition to determining the number of required scenarios that can be obtained as a corollary of our results presented so far.

\fi

\ifshort
\else

\begin{figure}
	\centering
	\includegraphics[width=0.65\textwidth]{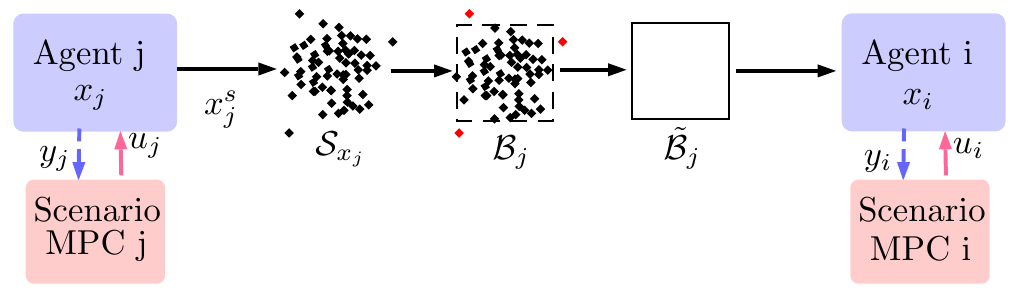}
	\caption{Pictorial representation of the proposed inter-agent soft communication scheme. 
		\myblue{$\cl{S}_{\bs{x}_j}$ is the set of ${S}_{i}$ scenarios, $\cl{B}_{j}$ is the parametrized set characterized in the optimization problem \eqref{privacy_box} using only $\tilde{S}_{j}$ scenarios, and  $\tilde{\cl{B}}_{j}$ is the solution of the optimization problem \eqref{privacy_box}. 
			The red dots refer to the difference between ${S}_{i}$ and $\tilde{S}_{j}$.} 
		The measurement variables $y_i$ and $y_j$ are the full state information $x_i$ and $x_j$, respectively, which are sent to the controllers.} 
	\label{fig:inter_agent}
	\vspace{-0.5cm}
\end{figure}

\fi

We now describe the soft communication protocol in more detail. 
The neighboring agent $j\in\cl{N}_i$ has to first generate $\tilde{S}_{j}$ samples of $\bs{x}_j$ in order to build the set $\cl{S}_{\bs{x}_j}$. 
It is important to notice that in the soft communication protocol the number $\tilde{S}_{j}$ of samples generated by agent $j$ may be different than the one needed by agent $i$, which is ${S}_{i}$, as will be remarked later.
Let us then introduce ${\cl{B}}_{j} \subset \R^{n_j} $ as a bounded set containing all the elements of $\cl{S}_{\bs{x}_j}$. 
We assume for simplicity that ${\cl{B}}_j$ is an axis-aligned hyper-rectangular set \cite{margellos2014road}. 
This is not a restrictive assumption and any convex set, e.g., ellipsoids and polytopes, could have been chosen instead as described in \cite{rostampour2017probabilistic}.
We can define ${\cl{B}}_{j} := [-\bs{b}_j, \bs{b}_j]$ as an interval, where the vector $\bs{b}_j\in\R^{n_j}$ defines the hyper-rectangle bounds.

Consider now the following optimization problem that aims to determine the set $\cl{B}_{j}$ with minimal volume:
\begin{align}\label{privacy_box}
\begin{cases}
\ \min\limits_{\bs{b}_j\in\R^{n_j}} & \| \bs{b}_j \|_1 \\
\ \ \ \text{s.t.} & \bs{x}_j^{l} \in [-\bs{b}_j, \bs{b}_j] \ , \ \forall \bs{x}_j^{l} \in\cl{S}_{\bs{x}_j} \\
& \  l = 1, \cdots, \tilde{S}_{j}
\end{cases} \ ,
\end{align}
where $\tilde{S}_{j}$ is the number of samples $\bs{x}_j \in\cl{S}_{\bs{x}_j}$ that neighboring agent $j$ has to take into account in order to determine  ${\cl{B}}_{j}$.
If we denote by $\tilde{\cl{B}}_j = [-\tilde{\bs{b}}_j, \tilde{\bs{b}}_j]$ the optimal solution of \eqref{privacy_box} computed by the neighbor agent $j$, then for implementing the soft communication protocol, agent $j$ needs to communicate only the vector $\tilde{\bs{b}}_j$ along with the level of reliability $\tilde{\alpha}_j$ to the agent $i$.

\begin{definition}\label{soft_box}
	A set $\tilde{\cl{B}}_j$ is called $\tilde{\alpha}_j-$reliable if 
	\begin{align}
	\Pb \left[\,\bs{x}_j \in \mathbb{X}_j \, : \, \bs{x}_j \notin [-\tilde{\bs{b}}_j, \tilde{\bs{b}}_j] \right] \leq 1-\tilde{\alpha}_j \,,
	\end{align}
	and we refer to $\tilde{\alpha}_j$ as the level of reliability of the set $\tilde{\cl{B}}_j$.
\end{definition}

We now provide the following theorem to determine $\tilde{\alpha}_j$ as the level of reliability of the set $\tilde{\cl{B}}_j$.
%
%
\begin{theorem}\label{reliable_thm}
	Fix $\tilde{\beta}_j\in(0,1)$ and let 
	\begin{align}\label{reliable_bound}
	\tilde{\alpha}_j = \sqrt[\kappa]{\frac{\tilde{\beta}_j}{ \binom{{\tilde{S}_{j}}}{n_j} } }\ .
	\end{align}
	\myblue{where $\kappa = {\tilde{S}_{j}} - n_j$ is the degree of the root}. We then have 
	\begin{align}\label{app_scenario}
	\Pb^{\tilde{S}_{j}} \Big\{ &\{ \bs{x}_j^1, \cdots, \bs{x}_j^{\tilde{S}_{j}} \} \in \mathbb{X}_j^{\tilde{S}_{j}} : \Pb \left[\,\bs{x}_j \in \mathbb{X}_j \, : \, \bs{x}_j \notin [-\tilde{\bs{b}}_j, \tilde{\bs{b}}_j] \right] \leq 1-\tilde{\alpha}_j 
	\Big\} \geq 1-\tilde{\beta}_j \,. 
	\end{align}
\end{theorem}

The proof is provided in the Appendix. \QEDB

Theorem \ref{reliable_thm} implies that given an hypothetical new sample $\bs{x}_j \in \mathbb{X}_j$, agent $j\in\cl{N}_i$ has a confidence of at least $1-\tilde{\beta}_j$ that the probability of $\bs{x}_j \in \tilde{\cl{B}}_j = [-\tilde{\bs{b}}_j, \tilde{\bs{b}}_j] $ is at least $\tilde{\alpha}_j$. 
Therefore, one can {\em rely on $\tilde{\cl{B}}_j$ up to $\tilde{\alpha}_j$ probability}.
The number of samples $\tilde{S}_{j}$ in the proposed formulation \eqref{privacy_box} is a design parameter chosen by the neighboring agent $j\in\cl{N}_i$. 
\ifshort
\else
In Fig.~\ref{fig:inter_agent} the number of red dots refers to the difference between ${S}_{i}$ and $\tilde{S}_{j}$.
\fi
We however remark that one can also set a given $\tilde{\alpha}_j$ as the desired level of reliability and obtain from \eqref{reliable_bound} the required number of samples $\tilde{S}_{i}$.

When an agent $i\in\cl{N}$ and its neighbor $j\in\cl{N}_i$ adopt the soft communication scheme, there is an important effect on the probabilistic feasibility of agent $i$, following \Cref{ccp_sub}.
Such a scheme introduces some level of stochasticity on the probabilistic feasibility of agent $i$, due to the fact that the neighboring information is only {\em probabilistically reliable}. 
This will affect the local probabilistic robustness guarantee of feasibility as it was discussed in \Cref{scenario2_thm} and consequently in \Cref{scenario_thm}.
To accommodate the level of reliability of neighboring information,  we need to marginalize a joint cumulative distribution function (cdf) probability of $\bs{x}_i$ and the generic sample $\bs{x}_j\in \mathbb{X}_j$ appearing in \Cref{reliable_thm}.
We thus have the following theorem, which can be regarded as the main theoretical result of this section.

\begin{theorem}\label{thm:probabilistically_reliable_threshold}
	Given $\tilde{\alpha}_j\in(0,1)$ for all $j\in\mathcal{N}_i$ and a fixed $\alpha_i=1-\varepsilon_i\in(0,1)$, the state trajectory of a generic agent \myblue{$i\in\cl{N}$} is probabilistically $\bar{\alpha}_i$--feasible
	for all $\bs{w}^{}_i \in {\cl{W}}_{i} \,,\, \bs{\delta}_{i} \in {\Delta}_{i}$,  i.e., 
	\begin{align}\label{new_chance_cons}
	\Pb\left[x_{i,k+\ell} \in \cl{X}_i \,,\, \ell \in \cl{T}_+ \right] \geq \bar{\alpha}_i \,,
	\end{align}
	where $\bar{\alpha}_i = 1 - \frac{1-\alpha_i}{\tilde{\alpha}_i}$ such that $\tilde{\alpha}_i=\prod_{j\in\cl{N}_i}(\tilde{\alpha}_j)$.
\end{theorem}

The proof is provided in the Appendix. \QEDB

Following the statement of \Cref{thm:probabilistically_reliable_threshold}, it is straightforward to observe that if for all neighboring agents $j\in\cl{N_i}$, $\tilde{\alpha}_j \rightarrow 1$ then $\bar{\alpha}_i \rightarrow \alpha_i$. This means that if {\em the level of reliability} of the neighboring information is one, i.e. $\Pb\big[\, \bs{x}_j\in\tilde{\cl{B}}_j \,  : \, \forall j\in\cl{N_i} \big] = 1$, then, the state feasibility of agent $i$ will have the same probabilistic level of robustness as the hard communication scheme, i.e. $\Pb\big[\, \bs{x}_i \in \mathbb{X}_i \, \big] \geq \alpha_i =1 -\varepsilon_i\,$. 
Combining this result with the statement of \Cref{scenario2_thm}, the proposed soft communication scheme introduces some level of stochasticity on the feasibility of the large-scale system as in \eqref{ccp_ls}.
In particular, $\varepsilon_i\in(0,1)$ the level of constraint violation in each agent $i\in\cl{N}$ will increase, since it is proportional with the inverse of $\prod_{j\in\cl{N}_i}(\tilde{\alpha}_j)\in(0,1)$, and therefore, $\varepsilon=\sum_{i\in\cl{N}}\varepsilon_i\in(0,1)$ will also increase.
After receiving the parametrization of $\tilde{\cl{B}}_j$ and the level of reliability  $\tilde{\alpha}_j$, agent $i\in\mathcal{N}$ should immunize itself against all possible variation of $\bs{x}_j\in\tilde{\cl{B}}_j$ by taking the worst-case of $\tilde{\cl{B}}_j$, similar to the worst-case reformulation proposed in \cite[Proposition 1]{rostampour2016robust}. 
It is important to notice that in this way, we decoupled the sample generation of agent $j\in\cl{N}_i$ from agent $i\in\cl{N}$.

\ifshort
\else

We summarize this section by mentioning that Fig.~\ref{fig:inter_agent} depicts a conceptual representation of the proposed soft communication scheme between two neighboring agents.
In the next section, we provide an operational framework that uses our developments in preceding sections in a more practical framework namely PnP distributed scenario MPC operations.

\fi

\myblue{
\section{Plug-and-Play Operational Framework}\label{PnP}

\ifshortt
\else

In this section, we develop a distributed framework to solve the centralized scenario MPC problem in \eqref{scp_ccp_ls} by using the so-called soft communication scheme developed in the preceding section. 

\fi

Using the soft communication scheme, each agent $i\in\mathcal{N}$ should first receive the parametrization of $\tilde{\cl{B}}_j$ from all its neighboring agents $j\in\mathcal{N}_i$ together with  the level of reliability  $\tilde{\alpha}_j$. 
Then, all agents should immunize themselves against all possible variation of $\bs{x}_j\in\tilde{\cl{B}}_j$.
To this end, we formulate the following \textit{robust-communication} scenario program for each agent $i\in\mathcal{N}$:
\begin{subequations}
	\label{scp_sub_robust}
	\begin{align}
	\min_{\bs{v}_i} \quad & J_i(\bs{x}_i, \bs{u}_i)  \\
	\text{s.t.} \quad &  {x}_{i,k+1}^{s} = A_{ii}(\delta_{i,k}^{s}) {x}_{i,k}^{s} + B_{i}(\delta_{i,k}^{s}) {u}_{i,k}^{s} + {q}_{i,k}^{s} \,,  \\
	&  {x}_{i,k+\ell}^{s} \in \cl{X}_i, \ell \in \cl{T}_+, \forall\bs{w}^{s}_i \in {\cl{S}}_{\bs{w}_i}, \forall\bs{\delta}_i^{s} \in {\cl{S}}_{\bs{\delta}_i}, \\
	&  {u}_{i,k}^{s} = K_{i} {x}_{i,k}^{s} + {v}_{i,k} \in \cl{U}_i \,,\, \forall k\in\cl{T} \,, \\
	&\begin{cases}\label{cons_robust}
	q_{i,k}^{s} = \sum\limits_{j\in\cl{N}_i} A_{ij}({\delta}_{i,k}^{s}) x_{j,k}^{s} + C_i({\delta}_{i,k}^{s}) w_{i,k}^{s}\\ 
	\forall\bs{x}_j^{s}\in\tilde{\cl{B}}_{j}\,, \forall j\in\cl{N}_i 
	\end{cases}.
	\end{align}
\end{subequations}
It is important to highlight that we refer to the aforementioned formulation as the \textit{robust-communication} scenario program, since the communicated variable between neighboring agents should be taken into consideration by the worst-case of $\tilde{\cl{B}}_j$. 
In this way, there is no need for the many iterations of \Cref{DSMPC}, and instead, a robustification against all possible variation of $\bs{x}_j\in\tilde{\cl{B}}_j$ is used.
It is also important to mention that another feature of using the soft communication scheme is the relaxation of the condition on the required number of scenarios, e.g., each agent $i\in\cl{N}$ requests $S_i$ number of scenarios from all its neighboring agents $j\in\cl{N}_i$, which may give rise to privacy concerns for the neighboring agents.  
The proposed problem \eqref{scp_sub_robust}, specifically the robust constraint \eqref{cons_robust}, can be solved using a robust (worst-case) reformulation similar to \cite[Proposition 1]{rostampour2016robust} and \cite{riverso2013plug}. 
It is important to highlight that based on Assumption~\ref{ass_scp_ls_feasible}, the proposed optimization problem \eqref{scp_sub_robust} should be feasible and  its feasibility domain has to be a nonempty domain.
In case of infeasible solution, each agent $i\in\mathcal{N}$ needs to generate new set of local scenarios, ${\cl{S}}_{\bs{w}_i}, {\cl{S}}_{\bs{\delta}_i}$, and also needs to receive a new set $\bs{x}_j\in\tilde{\cl{B}}_j$ from all its neighbors $j\in\mathcal{N}_i$.

We summarize the proposed distributed scenario MPC using soft communication scheme in \Cref{PnPSMPC}.
Similar to \Cref{DSMPC}, in \Cref{PnPSMPC} it is also assumed that the feedback control gain matrices $K_i$ for all agents $i\in\cl{N}$ are given \eqref{feedback_law}, and the coupling terms $A_{ij}(\delta_{i,k})$ are known between each agent $i\in\cl{N}$ and its neighboring agents $j\in\cl{N}_i$.
It is important to note that Step 5 of \Cref{PnPSMPC}, initializes $\tilde{\cl{B}}_{j}$ of all neighboring agents $j\in\cl{N_i}$, which is only used for the initial iteration in Step 8, and then, at each iteration all agents $i\in\cl{N}$ will send and receive $\tilde{\cl{B}}_{j}$ from all its neighboring agents $j\in\cl{N_i}$ as in Steps 11 and 12, respectively.

We also summarize the steps that are needed for plug-in and plug-out operations of each agent $i\in\cl{N}$ in \Cref{PnPframework}.
Note that in a plugged-in or plugged-out operation all agents $i\in\cl{N}$ have to update their $\varepsilon_i$ with $\beta_i$ to respect the condition in \Cref{scenario2_thm} as in \eqref{thm2condition} to achieve the desired level of constraint feasibility for the large-scale system \eqref{dyn_LS}.
One can also redesign $K_i$ to potentially improve the local control performance of each agent $i\in\cl{N}$.

\begin{algorithm}[t]
	\caption{\myblue{DSMPC using Soft Communication Scheme}}\label{PnPSMPC}
	\small
	\begin{algorithmic}[1]
		\State \textbf{Decompose} the large-scale dynamical system \eqref{dyn_LS} into $N$ agents as the proposed form in \eqref{dyn_sub}
		\State \textbf{Determine} the index set of neighboring agents $\cl{N_i}$ for each agent $i\in\cl{N}$
		\ForAll{$i\in\cl{N}$} \label{pnp}
		\State \textbf{fix} initial state ${x}_{i,0} \in \cl{X}_i$, $\varepsilon_i\in(0,1)$, and  $\beta_i\in(0,1)$ 
		\Statex \quad\,\ such that $\varepsilon = \sum\nolimits_{i\in\cl{N}}\varepsilon_i\in(0,1) \,,\, \beta = \sum\nolimits_{i\in\cl{N}}\beta_i\in(0,1)$ \label{thm2condition}
		\State \textbf{initialize} $\tilde{\cl{B}}_{j}$ for all neighboring agents $j\in\cl{N_i}$ 
		\State \textbf{determine} $\bar{S}_i\in(0,+\infty)$ to approximate the objective 
		\Statex \quad\,\ function, and $S_{i}$ following \Cref{scenario2_thm} to approximate 
		\Statex \quad\,\ the chance constraint \eqref{ccp_cons_sub} in an equivalent sense
		\State \textbf{generate} $\bar{S}_i$, $S_{i}$ scenarios of $\bs{w}_i$, $\bs{\delta}_i$ to determine the sets 
		\Statex \quad\,\ of $\bar{\cl{S}}_{\bs{w}_i,\bs{\delta}_i}$ and $\cl{S}_{\bs{w}_i}$, $\cl{S}_{\bs{\delta}_i}$ \label{start}
		\State  {\textbf{solve} the proposed optimization problem in \eqref{scp_sub_robust} by taking 
			\Statex \quad\,\ into account the worst-case of $\tilde{\cl{B}}_{j}$ and determine  $\bs{v}^*_i$}
		\State \textbf{generate} $\tilde{S}_{i}$ scenarios of $\bs{x}_i$ using the dynamical system 
		\Statex \quad\,\  of agent $i$ in form of \eqref{dyn_sub} and $\bs{v}^*_i$ together with $\cl{S}_{\bs{w}_i}$, $\cl{S}_{\bs{\delta}_i}$
		\State \textbf{determine} set $\tilde{\cl{B}}_{i}$ by solving the optimization problem \eqref{privacy_box} 
		\State \textbf{send} the set $\tilde{\cl{B}}_{i}$ to all neighboring agents $j\in\cl{N_i}$
		\State \textbf{receive} the sets $\tilde{\cl{B}}_{j}$ from all neighboring agents $j\in\cl{N_i}$
		\State \textbf{apply} the first input of solution $u_{i,k}^* = K_{i} x_{i,k} + v^*_{i,k}$ into 
		\Statex \quad\,\ the uncertain subsystem \eqref{dyn_sub} 
		\State \textbf{measure} the state and substitute it as the initial state of 
		\Statex \quad\,\ the next step ${x}_{i,0}$ 
		\State  \textbf{set} \ $k \gets k+1$ 
		\State  \textbf{go to} Step \eqref{start}
		\EndFor
	\end{algorithmic}
\end{algorithm}



\begin{algorithm}[t]
	\caption{Plug-and-Play Operation}\label{PnPframework}
	\begin{algorithmic}[1]
		\Statex \textbf{Plug-in Operation}
		\State \textbf{Add} the number of new subsystems into the previous number of subsystems, e.g. one additional agent $N \gets N+1$ such that $|\cl{N}| = N+1$
		\State \textbf{Update} the index set of neighboring agents $\cl{N_i}$ for each agent $i\in\cl{N}$
		\State \textbf{Go to} Step \eqref{pnp} of \Cref{PnPSMPC}
		\Statex \textbf{Plug-out Operation}
		\State \textbf{Remove} the number of departing subsystems from the previous number of subsystems, e.g. one agent leaving leads to $N \gets N-1$ such that $|\cl{N}| = N-1$
		\State \textbf{Update} the index set of neighboring agents $\cl{N_i}$ for each agent $i\in\cl{N}$
		\State \textbf{Go to} Step \eqref{pnp} of \Cref{PnPSMPC}
	\end{algorithmic}
\end{algorithm}

\subsection*{Comparison: DSMPC without vs. with Soft Communication}

We now provide a comparison between the proposed distributed scenario MPC without vs. with soft communication in \Cref{DSMPC} and \Cref{PnPSMPC}, respectively, in terms of their computational complexities and the conservatism level of obtained solutions as follows:

1) Computational complexity: 
It is easy to realize that the proposed distributed scenario MPC (\Cref{PnPSMPC}) using soft communication scheme is computationally advantageous compared to the distributed scenario MPC without soft communication (\Cref{DSMPC}). 
This can be clearly seen by just comparing Step 7 and Step 8 in \Cref{DSMPC} and \Cref{PnPSMPC}, respectively.  
In Step 7 of \Cref{DSMPC} each agent $i\in\mathcal{N}$ needs to execute the proposed distributed scenario exchange scheme in \Cref{alg:distributed scenario exchange ADMM} until all its neighboring agents agree on consistency of the set of their state variable scenarios $\cl{S}_{\bs{x}_j}$.
Whereas in Step 8 of \Cref{PnPSMPC} each agent $i\in\mathcal{N}$ needs to solve only once the robust-communication scenario program in \eqref{scp_sub_robust} and also solve the proposed optimization problem \eqref{privacy_box} in Step 10 of \Cref{PnPSMPC} in order to determine the set $\tilde{\cl{B}}_{i}$.

2) Conservatism level: 
As it is shown in \Cref{scenario2_thm}, the proposed distributed scenario MPC using \Cref{DSMPC} retrieves exactly the same property of the original centralized MPC \eqref{scp_ccp_ls} under certain conditions, i.e., $\varepsilon = \sum_{i\in\cl{N}} \varepsilon_i \in (0,1)$, $\beta = \sum_{i\in\cl{N}} \beta_i \in (0,1)$. 
This is however different in the distributed scenario MPC using the proposed soft communication scheme as it is presented in \Cref{thm:probabilistically_reliable_threshold}.
In fact such a scheme introduces some level of stochasticity on the probabilistic feasibility of the agent $i\in\mathcal{N}$, due to the probabilistic reliability of the neighboring information. 
It is important to mention that guaranteeing the optimality of the obtained solutions in terms of performance objective(s) using the proposed distributed scenario MPC without or with soft communication in \Cref{DSMPC} and \Cref{PnPSMPC}, respectively, is not included in the scope of this paper and it is subject of our ongoing research work.

}

\section{Numerical Study}\label{sim}

\myblue{
	We simulate four different problem formulations, namely a centralized SMPC (CSMPC) using \eqref{scp_ls}, a distributed SMPC (DSMPC) via \Cref{DSMPC}, and DSMPC with the proposed soft communication scheme with $0.85-$reliability (DSMPCS$-0.85$) as described in \Cref{soft_box} and DSMPCS$-0.50$, both following \Cref{PnPSMPC} in a closed-loop control system framework.
	For comparison purposes, we also present the results obtained via decoupled SMPC (DeSMPC), where the impact of coupling neighboring dynamics in \eqref{dyn_sub} are relaxed.
}

In Fig.~\ref{fig:coupled_viol_rooms}, \ifshortt\else Fig.~\ref{fig:coupled_viol}\fi and Fig.~\ref{fig:zoomin_coupled_viol}, we evaluate our proposed framework in terms of a-posteriori feasibility validation of the obtained results in both case studies. 
The "red" line represents the results obtained via DeSMPC, the "blue" line shows the results obtained via CSMPC, the "magenta" presents the results obtained by using DSMPC, the "dark green" and "light green" lines show the results obtained via DSMPCS$-0.85$ and DSMPCS$-0.50$, respectively.
The "black" lines indicate the bounds of the three dynamically coupled systems.

\subsection{Three-Room Case Study}\label{three_room}

We simulate a building climate comfort system with three rooms \ifshort\else shown in Fig.~\ref{fig:ThreeRoom} \fi such that the temperature of rooms are dynamically coupled without any common constraints.
The outside weather temperature and the related looses, e.g. through windows, is considered as the private uncertainty.

Consider now a three-room building system dynamics: ${x}_{k+1} = A(\delta_k) {x}_{k} + B(\delta_k) {u}_{k} + C(\delta_k) w_{k}$\,, where
$B = \text{diag}([0.01, 0.01, 0.01])$, $C = \text{diag}([0.02, 0.02, 0.02])$, and $A$ as follows:
\begin{equation*}
\begin{aligned}
A = \begin{bmatrix}
0.2 & 0.3 & 0 \\
0.2 & 0.1 & 0.1 \\ 
0.2 & 0 & 0.4
\end{bmatrix}
\end{aligned}
\end{equation*}
such that $A(\delta_k) = A + \delta_k$ and $B(\delta_k) = B + \delta_k$  as well as $C(\delta_k) =  C + \delta_k$.
The system matrices are a simplified model of a three-room building such that the states $x_{i,k}$ for $i=1,2,3$, denote the temperature of rooms.
The uncertain variable $\delta_k\in\R$ represents the modeling errors, losses through windows, and $w_k\in\R$ can be realized as the outside weather temperatures.

To generate random scenarios from the additive disturbance, we built a discrete normal process such that one day hourly-based forecasted (nominal) outside weather temperature is used which varies within $10\%$ of its nominal scenario at each sampling time.
As for the uncertainty $\delta_k$, we generate a random variable from a normal distribution with a mean value $0$, variance $1$ and a maximal magnitude of $0.01$ at each sampling time.

The initial state variables are $[21 \ 19 \ 23]^\top$ and the objective is to keep the temperature of rooms within our desired lower $[20.5 \ 18.5 \ 22.5]^\top$ and $[21.5 \ 19.5 \ 23.5]^\top$ upper bounds at the minimum control input $u_k$.
The control input $u_k$ is also constrained to be within $-1.5$ [kWh] and $1.5$ [kWh] for all three rooms, due to actuator saturation.
The initialization of the $\tilde{\cl{B}}_{j}$ for all neighboring agents $j\in\cl{N_i}$ as in Step 5 in \Cref{PnPSMPC} can be done for instance by assuming the initial temperature of the neighboring rooms are known for all rooms.

\ifshort
\else

\begin{figure}
	\centering
	\includegraphics[width=0.5\textwidth]{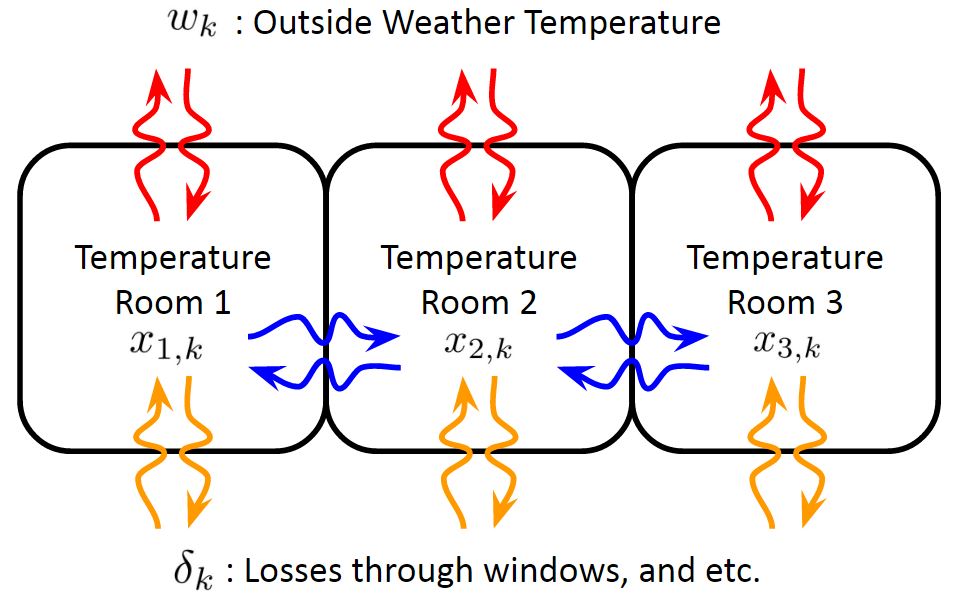}
	\caption{An example of three-room building climate comfort system. $x_{1,k}$, $x_{2,k}$ and $x_{3,k}$ are states of the building that are related to the temperature dynamics of rooms. The temperature dynamics are influenced by the outside weather temperature $(w_k)$ shown by `red' arrows and the possible losses $(\delta_k)$ through windows, and etc., of the rooms which are represented via `orange' arrows. Moreover, the `blue' arrows denotes the impact of the neighboring rooms on each other.}
	\label{fig:ThreeRoom}
\end{figure}

\fi

\begin{figure}
	\centering
	\includegraphics[width=0.65\textwidth]{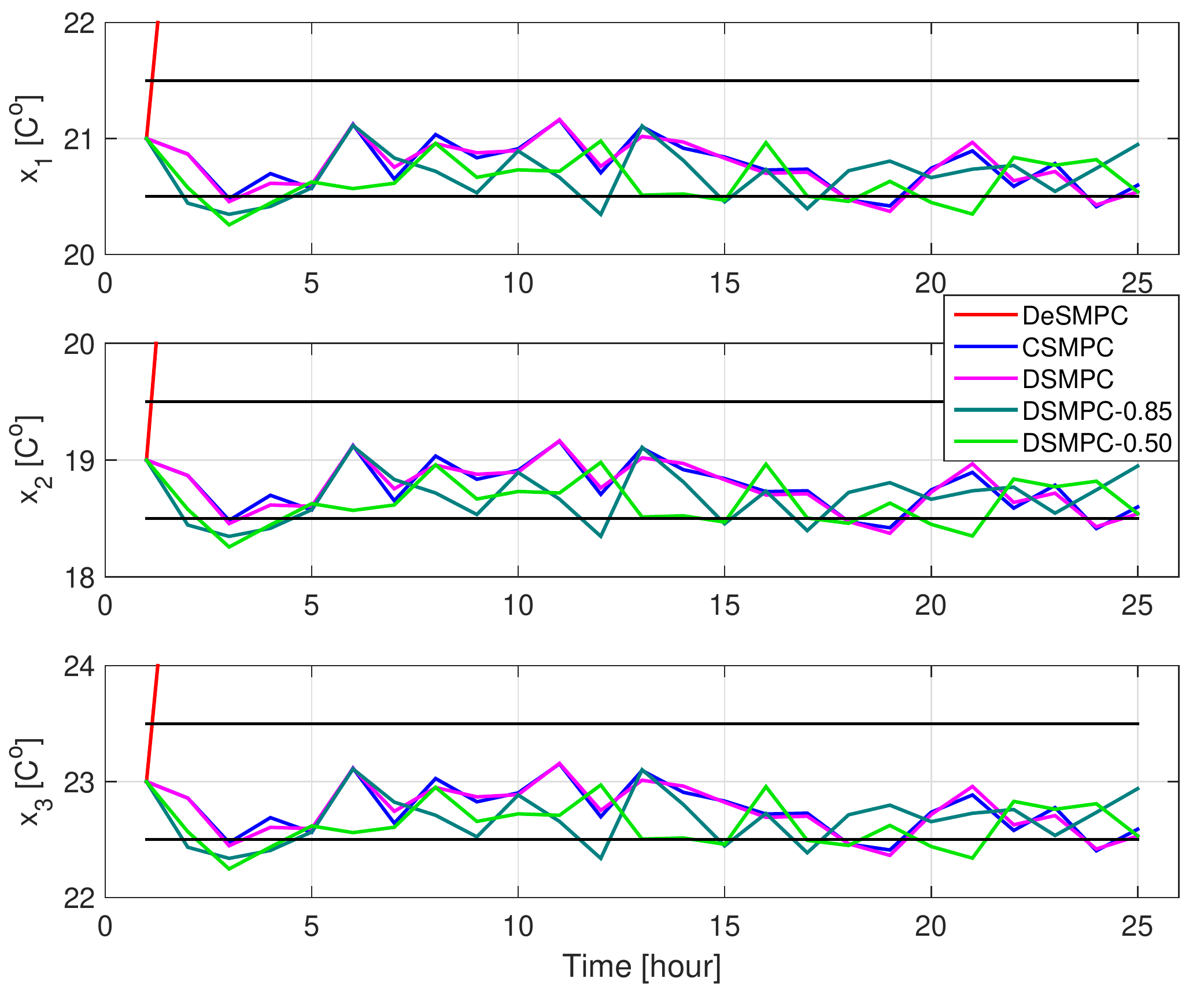}
	\caption{A-posteriori feasibility validation of the obtained results. The "red", "blue", "magenta", "dark green", and "light green" lines are related to the obtained results via DeSMPC, CSMPC, DSMPC, DSMPCS$-0.85$, and DSMPCS$-0.50$, respectively. The "black" lines are related to the upper bound values.
		The top, middle, and down plots are related to the temperatures of room $1, 2$, and $3$, respectively.}	
	\label{fig:coupled_viol_rooms}
\end{figure}

In Fig.~\ref{fig:coupled_viol_rooms}, the dynamically coupled state trajectories for all three rooms are shown.
One can clearly see in Fig.~\ref{fig:coupled_viol_rooms} that the dynamically coupled state trajectories are feasible in a probabilistic sense, since the agent operations are within the lower and upper bounds compared to DeSMPC that violates the constraints completely; the obtained solution via DeSMPC is completely outside of the feasible areas after the first sampling time and thus, we just keep the other results for our discussions.
This is a direct result of \Cref{scenario2_thm} such that the obtained solutions via our proposed formulations have to be probabilistically feasible, that can be clearly seen  in Fig.~\ref{fig:coupled_viol_rooms}, since the trajectories are on the lower bounds.

\subsection{Three-Building (ATES Systems) Case Study}\label{ATES_STGs}

We next simulate the thermal energy demands of three buildings modeled using realistic parameters and the actual registered weather data in the city center of Utrecht, The Netherlands, where these buildings are located and these had been equipped with aquifer thermal energy storage (ATES) systems.  
\ifshortt
A simulation study is carried out for one year based on actual weather conditions from March 2011-March 2012 in a weekly-based time steps with three months prediction horizon to meet the coupling constraints between ATES systems.
We refer interested readers to \cite{rostampour2017probabilistic} for the detailed explanations on this case study.

\else
An ATES system consists of two wells (warm and cold water wells) and it is considered as a heat source or sink, or as a storage for thermal energy that operates in a seasonal mode.
A large-scale network of interconnected buildings, that are constrained via the state variables of ATES systems, which are the volume of water and the thermal energy content of each well, was modeled in \cite{rostampour2017probabilistic}.
To prevent overlap of nearby systems there are constraints on a growing thermal radius, $r_{i,k}^{h} \, \text{[m]}$, $r_{i,k}^{c} \, \text{[m]}$, of each well of an ATES system with its neighboring agents, e.g. $r_{i,k}^{h} + r_{j,k}^{c} \leq d_{ij}$ for each agent $i\in\cl{N}$ and its neighboring agent $j\in\cl{N}_i$.
In \cite[Corollary 1]{rostampour2017probabilistic}, it was shown that such a constraint on the growing thermal radius can be replaced with a constraint on the state variables (volume of water) of the ATES systems. 
We refer interested readers to \cite{rostampour2017probabilistic} for the detailed explanations on this case study.

A simulation study is carried out for one year based on actual weather conditions from March 2011-March 2012 in a weekly-based time steps with three months prediction horizon to meet the coupling constraints between ATES systems.
We introduce additive disturbance and uncertainty sources into the deterministic first-order dynamical model of \cite{rostampour2017probabilistic}.
It has been reported in \cite{rostampour2016control} that the ambient temperature of water in aquifers is changing over the long life-time usage of ATES systems.
We capture such a behavior by introducing an additive unknown (disturbance) source of energy which yields a time correlation function via the initial value of energy content variable of an ATES system. 
In addition to this, the system parameters of an ATES system are a function of the stored water temperature in each well, e.g. see \cite[Fig.~2]{rostampour2017probabilistic}. 
We therefore introduce a small level of perturbation as an uncertainty in the parameters of the ATES system dynamics. 

\fi

To generate random scenarios from the additive disturbance, we built a discrete normal process such that the nominal scenario is $10\%$ of the amplitude of the energy content in a deterministic ATES system model and varies within $5\%$ of its nominal scenario at each sampling time.
As for the uncertainty $\delta_k$, we generate a random variable from a Gaussian distribution with a mean value $0$, variance $0.3$ and a maximal magnitude of $0.03$ at each sampling time.

\ifshortt
\else

\begin{figure}
	\centering
	\includegraphics[width=0.65\textwidth]{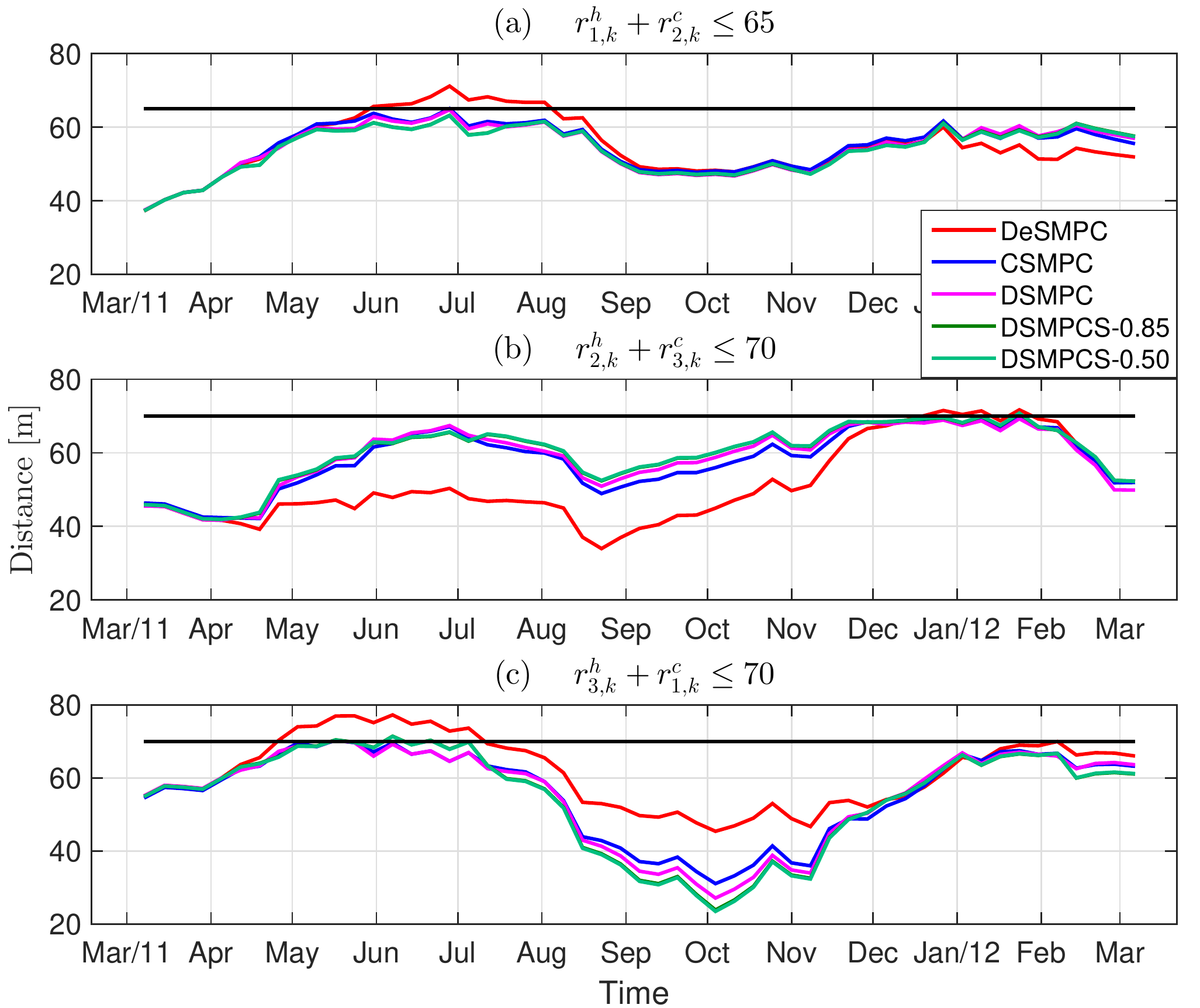}
	\caption{A-posteriori feasibility validation of the obtained results. The "red", "blue", "magenta", "solid green", and "dashed green" lines are related to the obtained results via DeSMPC, CSMPC, DSMPC, DSMPCS$-0.85$, and DSMPCS$-0.50$, respectively. The "black" lines are related to the upper bound values.}
	\label{fig:coupled_viol}
\end{figure}

\fi

\begin{figure}
	\centering
	\includegraphics[width=0.65\textwidth]{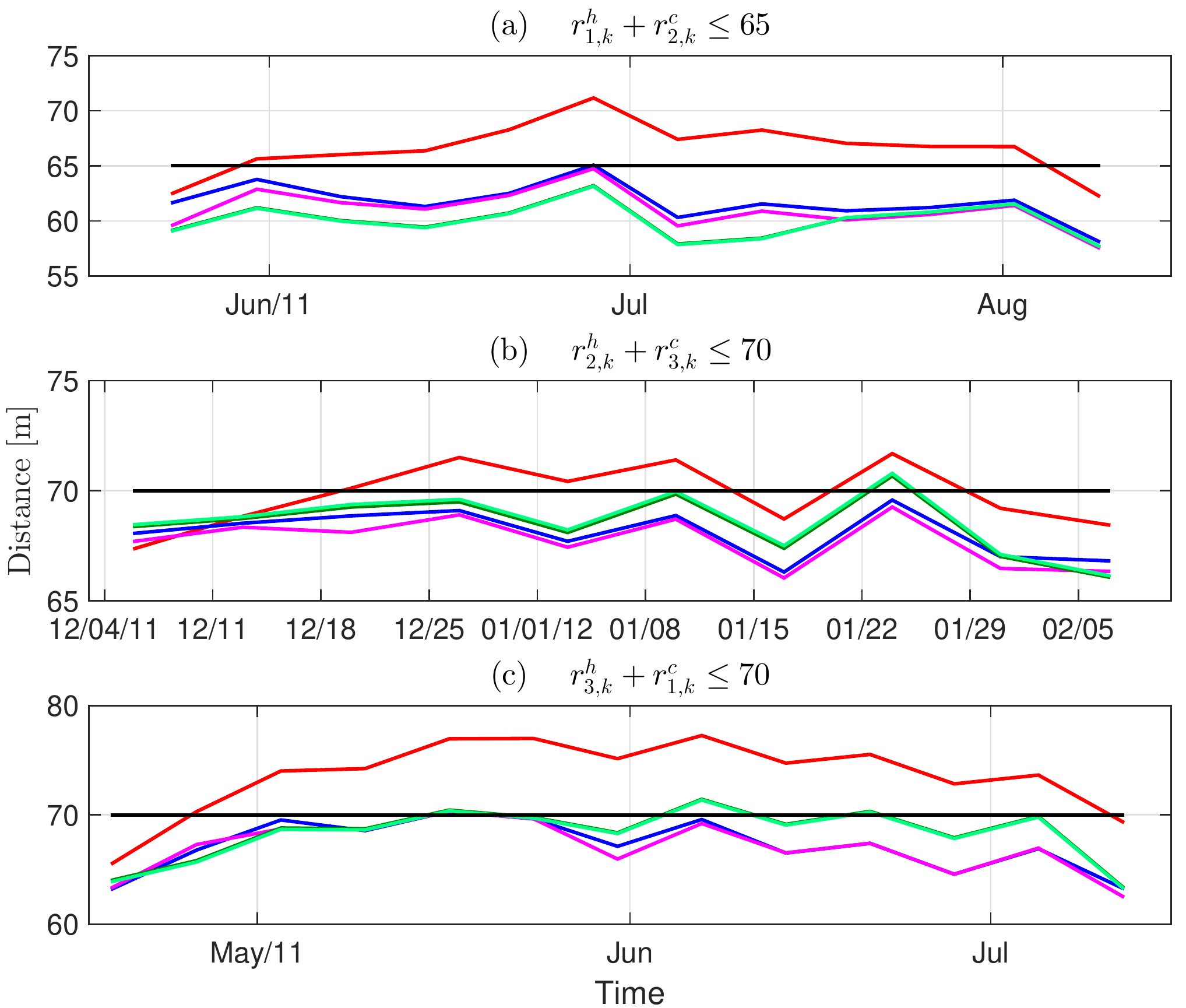}
	\ifshortt
	\caption{A-posteriori feasibility validation of the obtained results. The "red", "blue", "magenta", "solid green", and "dashed green" lines are related to the obtained results via DeSMPC, CSMPC, DSMPC, DSMPCS$-0.85$, and DSMPCS$-0.50$, respectively. The "black" lines are related to the upper bound values.}
	\else
	\caption{Zoom-in of the critical time periods in \Cref{fig:coupled_viol}.}
	\fi
	\label{fig:zoomin_coupled_viol}
\end{figure}

\ifshortt

In Fig.~\ref{fig:zoomin_coupled_viol} we show a-posteriori feasibility validation of the coupling constraints between each agent $i = 1, 2, 3$, and neighboring agents, such that (a) shows the results from mid-May to mid-August 2011, (b) presents the results of December 2011 to February 2012, and (c) depicts the results of mid-April to mid-July 2011. 
As a first desired achievement, one can clearly see in Fig.~\ref{fig:zoomin_coupled_viol} that the constraints are feasible in a probabilistic sense, since the agent operations are quite close to the upper bounds in the critical time periods compared to DeSMPC that violates the constraints. 
Strictly speaking, using our proposed framework one can achieve the maximum usage of the aquifer (subsurface) to store thermal energy without affecting the neighboring thermal storage.  
This is a direct result of \Cref{scenario2_thm} such that the obtained solutions via our proposed formulations have to be probabilistically feasible.

\else

In Fig.~\ref{fig:coupled_viol} we show a-posteriori feasibility validation of the coupling constraints between each agent $i = 1, 2, 3$, and neighboring agents, e.g. $r_{1,k}^h + r_{2,k}^c \leq 65$, $r_{2,k}^h + r_{3,k}^c \leq 70$, and $r_{3,k}^h + r_{1,k}^c \leq 70$.
Fig.~\ref{fig:zoomin_coupled_viol} focuses on the critical time periods in Fig.~\ref{fig:coupled_viol}, where neighboring agents are injecting thermal energies with different pump flow rates.  
Fig.~\ref{fig:zoomin_coupled_viol}(a) shows the results from mid-May to mid-August 2011, Fig.~\ref{fig:zoomin_coupled_viol}(b) presents the results of December 2011 to February 2012, and Fig.~\ref{fig:zoomin_coupled_viol}(c) depicts the results of mid-April to mid-July 2011. 

\fi

\ifshortt
\else

As a first desired achievement, one can clearly see in Fig.~\ref{fig:coupled_viol} that the constraints are feasible in a probabilistic sense, since the agent operations are quite close to the upper bounds in the critical time periods compared to DeSMPC that violates the constraints. 
Strictly speaking, using our proposed framework one can achieve the maximum usage of the aquifer (subsurface) to store thermal energy without affecting the neighboring thermal storage.  
This is a direct result of \Cref{scenario2_thm} such that the obtained solutions via our proposed formulations have to be probabilistically feasible, that can be clearly seen  in Fig.~\ref{fig:zoomin_coupled_viol}, since the trajectories are on the upper bounds.

\fi

It is worth to mention that both Fig.~\ref{fig:coupled_viol_rooms} and Fig.~\ref{fig:zoomin_coupled_viol} illustrate our other two main contributions more precisely: 
1) the obtained results via CSMPC (blue line) and DSMPC (magenta line) are practically equivalent throughout the simulation studies; 
this is due to \Cref{propos:decomp_exact} and \Cref{scenario2_thm}. 
Actually, the solutions via DSMPC are slightly more conservative compared to the results via CSMPC, and this is a direct consequence of \Cref{scenario2_thm}. 
In fact the level of violation in CSMPC is considered to be $\varepsilon = 0.05$ and leading to $\varepsilon_i = 0.0167$ for all agents, which is more restrictive. 
2) the proposed soft communication scheme yields less conservative solutions as explicitly derived in \Cref{thm:probabilistically_reliable_threshold}, and can be clearly seen in Fig.~\ref{fig:coupled_viol_rooms} and Fig.~\ref{fig:zoomin_coupled_viol} with the obtained results via DSMPCS$-0.85$ (dark green) and DSMPCS$-0.50$ (light green).
Following \Cref{thm:probabilistically_reliable_threshold} the new violation level using DSMPCS$-0.85$ is $\bar{\varepsilon}_i = 0.0231$, and using DSMPCS$-0.50$ is $\bar{\varepsilon}_i = 0.0668$.
It is important to notice that the violation level of global chance constraint will increase to $\bar{\varepsilon} = 0.0702$ and $\bar{\varepsilon} = 0.2004$ using DSMPCS$-0.85$ and DSMPCS$-0.50$, respectively.

\section{Concluding Remarks and Future Work}\label{final}

In this paper, we presented a rigorous approach to distributed stochastic SMPC using the scenario-based approximation for large-scale linear systems with private and common uncertainty sources.
\myblue{
	We extended the existing results to quantify the robustness of the resulting solutions for both cases of private and common uncertainties in a distributed SMPC framework using the so-called distributed scenario exchange scheme based on ADMM.
	We then provided a novel inter-agent soft communication scheme to minimize the amount of information exchange between each subsystem.
	Using a set-based parametrization technique, we introduced a reliability notion and quantified the level of feasibility of the obtained solutions via the PnP distributed  SMPC integrated with the so-called soft communication scheme in a probabilistic sense.
}
The theoretical guarantees of the proposed distributed SMPC framework coincide with its centralized counterpart.

\ifshortt
\else

Our current research direction concentrates on enhancing the proposed framework to formally address the recursive feasibility and stability of the closed-loop.
It is important to point out that the discussion on scenario optimization was limited to constraint satisfaction probability.
\myblue{
	An interesting future work is to analyze the performance objective(s) by means of optimality of the scenario program solutions in a distributed setting, which is the subject of our ongoing research work.
}  
Finally, another interesting future research direction is to investigate the possibility of extending our proposed frameworks to the fault detection and isolation problems following the results in \cite{rostampour2017set}.

\fi

\section*{Appendix}

\noindent\textbf{Proof of \Cref{propos:decomp_exact}.}
Given Assumption~\ref{ass:decompose} and following the proposed decomposition technique in \Cref{dist_prob_state}, any optimizer of each subsystem $\bs{v}_i^*$ yields a feasible pair of the state and control input variables of its subsystem $\{\bs{x}_i^*, \bs{u}_i^*\}\in \mathbb{X}_i \times \mathbb{U}_i$ such that $\mathbb{X}_i = \cl{X}_i^T$, and $\mathbb{X}_i = \cl{U}_i^T$. 
Therefore, the collection of the optimizers $\bs{v}^* = \text{col}_{i\in\cl{N}}(\bs{v}_i^*)$ will yield the collection of feasible pairs of the state and control input variables of their subsystem:
\begin{align*}
\left\{\bs{x}^* = \text{col}_{i\in\cl{N}}(\bs{x}_i^*)\,, \bs{u}^* = \text{col}_{i\in\cl{N}}(\bs{u}_i^*)\right\}\in \mathbb{X} \times \mathbb{U} \,,
\end{align*}
where $\mathbb{X} := \cl{X}^T = \prod_{i\in\cl{N}}\mathbb{X}_i$ and $\mathbb{U} := \cl{U}^T = \prod_{i\in\cl{N}} \mathbb{U}_i$,
which eventually yields a feasible point for the optimization problem in \eqref{scp_ls}.
It is straightforward to use the above relation and to show that any optimizer of the optimization problem in \eqref{scp_ls} $\bs{v}^*$ also yields a feasible solution for the proposed optimization problem in \eqref{scp_sub}. 
We then have to show that both optimization problems will have the same performance index in terms of their objective function values.
Due to the proposed decomposition technique \myblue{and under Assumption~\ref{ass:decompose}}, it is easy to see that the objective function in \eqref{scp_ls} can be formulated as additive components such that each component represents the objective function of each subsystem $i\in\cl{N}$, and thus: 
$J(\bs{x}^*, \bs{u}^*) = \sum_{i\in\cl{N}} J_i(\bs{x}_i^*, \bs{u}_i^*) \,.$
The proof is completed.
\QEDB

\noindent\textbf{Proof of \Cref{scenario2_thm}.}
Define $\xi_{i,k} := (w_{i,k},\delta_{i,k},\{x_{j,k}\}_{j\in\mathcal{N}_i})$ to be a concatenated uncertain variable for each agent $i\in\mathcal{N}$ such that $\bs{\xi}_i := \{\xi_{i,k}\}_{k\in\mathcal{T}}$ is defined on probability space $(\Xi_i, \mathfrak{B}(\Xi_i),\Pb_{\bs{\xi_i}})$, where $\Pb_{\bs{\xi_i}}$ is a probability measure defined over $\Xi_i:= \mathcal{W}_i\times \Delta_i\times \prod_{j\in\mathcal{N}_i}\mathbb{X}_j$ and $\mathfrak{B}(\cdot)$ denotes a Borel $\sigma$-algebra.
Following this definition, it is straightforward to consider $\bs{\xi} = \text{col}_{i\in\mathcal{N}}(\bs{\xi}_i)$ and $\Xi = \prod_{i\in\mathcal{N}}\Xi_i$. 
Consider also the sample set $\mathcal{S}_{\bs{\xi}_i} = \mathcal{S}_{\bs{w}_i} \times \mathcal{S}_{\bs{\delta}_i}\times \prod_{j\in\mathcal{N}_i}\mathcal{S}_{\bs{x}_j}$ for each agent $i\in\mathcal{N}$ such that $\mathcal{S}_{\bs{\xi}} = \prod_{i\in\mathcal{N}} \mathcal{S}_{\bs{\xi}_i}$.

Consider now $\bs{v}^*$ to be the optimizer of the centralized scenario MPC problem \eqref{scp_ls} and define $\text{Vio}(\bs{v}^*)$ to be the violation probability of the chance constraint \eqref{ccp_ls} as follows:
\begin{subequations}
	\label{scenario_bound_statements}
	\begin{align}
	\text{Vio}(\bs{v}^*) := \mathbb{P}_{\bs{\xi}}\left[\, \bs{\xi}\in\Xi \,:\, g(\bs{v}^*,\bs{\xi}) \notin\mathbb{X} \,\right] \,, 
	\end{align}
	where $g(\bs{v}^*,\bs{\xi})$ represents the predicted state trajectory of large-scale system dynamics \eqref{dyn_LS} in a more compact form.
	In particular, the violation probability can be precisely written as $\text{Vio}(\bs{v}^*) := \Pb[\, \bs{w}\in\cl{W}, \bs{\delta}\in\Delta \,:\, {x}_{k+\ell} = A_{cl}(\delta_k) {x}_{k} + B(\delta_k) {v}_{k}^* + C(\delta_k) {w}_{k} \notin\cl{X} \,,\, \ell \in \cl{T}_+ \big\vert \, {x}_k = {{x}}_0 \,] \,,$ where $A_{cl}(\delta_k) = A(\delta_k) + B(\delta_k) K$.
	Following \Cref{scenario_thm}, we have
	\begin{align}
	\mathbb{P}^{S}_{\bs{\xi}}\left[\, \mathcal{S}_{\bs{\xi}}\in\Xi^{S}\,:\, \text{Vio}(\bs{v}^*) \leq \varepsilon \,\right] \geq 1-\beta \,.
	\end{align}
\end{subequations}
Given \Cref{propos:decomp_exact}, the problem \eqref{scp_sub} is an exact decomposition of the problem \eqref{scp_ls}.
This yields the following equivalence relations:
\begin{align}\label{decomp_relations}
\bs{v}^* := \text{col}_{i\in\mathcal{N}}(\bs{v}^*_i) \ , \ \mathbb{X} := \prod\nolimits_{i \in \mathcal{N}} \mathbb{X}_{i} \,,
\end{align}
where $\bs{v}^*$ is the optimizer of the problem \eqref{scp_ls} with $\bs{v}^*_i$ as the optimizer of each agent $i \in \mathcal{N}$ obtained via the problem \eqref{scp_sub}.
To this end, it is necessary to prove that the above statements \eqref{scenario_bound_statements} still hold under the aforementioned relations \eqref{decomp_relations}.
We now break down the proof into two steps.
We first show the results for each agent $i \in \mathcal{N}$, and then extend into the large-scale scenario MPC problem \eqref{scp_ls}. 

1) Define $\text{Vio}(\bs{v}^*_i)$ to be the violation probability of each agent $i \in \mathcal{N}$ fo the chance constraint \eqref{ccp_cons_sub} as follows:
\begin{align}
\text{Vio}(\bs{v}^*_i) := \mathbb{P}_{\bs{\xi}_i}\left[\, \bs{\xi}_i\in\Xi_i \,:\, g_i(\bs{v}^*_i,\bs{\xi}_i) \notin\mathbb{X}_i \,\right] \,, 
\end{align}
where $g_i(\bs{v}^*_i,\bs{\xi}_i)$ corresponds to the predicted state trajectory of subsystem dynamics \eqref{dyn_sub} for each agent $i \in \mathcal{N}$.
Applying the existing results in \Cref{scenario_thm} for each agent $i \in \mathcal{N}$, we have
\begin{align}
\mathbb{P}^{{S}_{i}}_{\bs{\xi}_i}\left[\, \mathcal{S}_{\bs{\xi}_i}\in\Xi_i^{{S}_{i}}\,:\, \text{Vio}(\bs{v}^*_i) \leq \varepsilon_i \,\right] \geq 1-\beta_i \,.
\end{align}

2) Following the relations \eqref{decomp_relations}, it is easy to rewrite $\text{Vio}(\bs{v}^*)$ in the following form:
\begin{align*}
\text{Vio}(\bs{v}^*) = \mathbb{P}_{\bs{\xi}}\left[\, \bs{\xi}\in\Xi \,:\, g(\bs{v}^*,\bs{\xi}) \notin\prod_{i\in\mathcal{N}}\mathbb{X}_i \,\right] \,. 
\end{align*}
It is then sufficient to show that for $S = \max_{i \in \mathcal{M}} S_{i}$:
\begin{align}\label{desired_bound}
\mathbb{P}^{{S}_{}}_{\bs{\xi}}\left[\, \mathcal{S}_{\bs{\xi}}\in\Xi^{{S}_{}}\,:\, \text{Vio}(\bs{v}^*) \geq \varepsilon \,\right] \leq \beta \,.
\end{align}
where $\varepsilon = \sum_{i\in\mathcal{N}}\varepsilon_i \in(0,1)$ and $\beta = \sum_{i\in\mathcal{N}}\beta_i \in(0,1)$. 
Hence
\begin{align*}
\text{Vio}(\bs{v}^*) &= \mathbb{P}_{\bs{\xi}}\left[\, \bs{\xi}\in\Xi \,:\, g(\bs{v}^*,\bs{\xi}) \notin\prod_{i\in\mathcal{N}}\mathbb{X}_i \,\right] \, \\
&= \mathbb{P}_{\bs{\xi}}\left[\, \bs{\xi}\in\Xi \,:\, \exists\, i\in\mathcal{N} \,,\, g(\bs{v}^*_i,\bs{\xi}_i) \notin\mathbb{X}_i \,\right] \, \\
&= \mathbb{P}_{\bs{\xi}}\left[\, \bigcup_{i \in\mathcal{N}}\left\{ \, \bs{\xi}_i\in\Xi_i \,:\,  g(\bs{v}^*_i,\bs{\xi}_i) \notin\mathbb{X}_i \, \right\} \right] \, \\
&\leq \sum_{i \in\mathcal{N}} \mathbb{P}_{\bs{\xi}_i}\left[\, \bs{\xi}_i\in\Xi_i :  g(\bs{v}^*_i,\bs{\xi}_i) \notin\mathbb{X}_i \,  \right] = \sum_{i \in\mathcal{N}} \text{Vio}(\bs{v}^*_i) . 
\end{align*}
The last statement implies that $\text{Vio}(\bs{v}^*) \leq \sum_{i \in\mathcal{N}} \text{Vio}(\bs{v}^*_i)$, and thus, we have
\begin{align*}
\mathbb{P}^{{S}_{}}_{\bs{\xi}}&\left[\, \mathcal{S}_{\bs{\xi}}\in\Xi^{{S}_{}}\, :\, \text{Vio}(\bs{v}^*) \geq \varepsilon \,\right]\\
&\leq \mathbb{P}^{{S}_{}}_{\bs{\xi}}\left[\, \mathcal{S}_{\bs{\xi}}\in\Xi^{{S}_{}}\,:\, \sum_{i \in\mathcal{N}} \text{Vio}(\bs{v}^*_i) \geq \sum_{i \in\mathcal{N}}\varepsilon_i \,\right] \\
&= \mathbb{P}^{{S}_{}}_{\bs{\xi}}\left[\,\bigcup_{i \in\mathcal{N}}  \left\{ \mathcal{S}_{\bs{\xi}_i}\in\Xi_i^{{S}_{i}}\,:\, \text{Vio}(\bs{v}^*_i) \geq \varepsilon_i \,\right\} \, \right] \\
&\leq \sum_{i \in\mathcal{N}} \mathbb{P}^{{S}_{i}}_{\bs{\xi}_i}\left[\, \mathcal{S}_{\bs{\xi}_i}\in\Xi_i^{{S}_{i}}\,:\, \text{Vio}(\bs{v}^*_i) \geq \varepsilon_i \, \right] \leq  \ \sum_{i \in\mathcal{N}}  \beta_i = \beta \ . 
\end{align*}
The obtained bounds in the above procedure are the desired assertions as it is stated in the theorem.
The proof is completed by noting that the feasible set  $\mathbb{X} = \prod\nolimits_{i \in \mathcal{N}} \mathbb{X}_{i}$ of the problem \eqref{scp_ls} and \eqref{scp_sub} has a {non-empty interior}, and it thus admits at least one feasible solution 
$\bs{v}^* = \text{col}_{i\in\mathcal{N}}(\bs{v}^*_i)$.
\QEDB

\noindent\textbf{Proof of \Cref{reliable_thm}.}
\Cref{app_scenario} is a direct result of the scenario approach theory in \cite{calafiore2006scenario}, if $\tilde{\beta}_j$ is chosen such that 
\begin{align}
{\binom{{\tilde{S}_{j}}}{n_j}} \tilde{\alpha}_j^{{\tilde{S}_{j}} - n_j} \leq \tilde{\beta}_j \,.
\end{align}
Considering the worst-case equality in the above relation and some algebraic manipulations, one can obtain the above assertion.
The proof is completed.
\QEDB

\noindent\textbf{Proof of \Cref{thm:probabilistically_reliable_threshold}.}
The proof consists of two main steps. 
We first provide an analytical expression for the robustness of the solution in agent $i$ by taking into account the effect of just one neighboring agent $j\in\cl{N}_j$, and then extend the obtained results for the case when the agent $i$ interacts with more neighboring agents, e.g. for all $j\in\cl{N}_j$.

Following \Cref{ccp_sub}, we have the following updated situation:
\begin{align*}
\alpha_i\leq \Pb\left[\, \bs{x}_i\in \mathbb{X}_i \ , \ \bs{x}_j\in  \tilde{\cl{B}}_j \, \right] \ ,
\end{align*}
which is a joint probability of $\bs{x}_i\in \mathbb{X}_i$ and $\bs{x}_j\in  \tilde{\cl{B}}_j$.
Such a joint probability can be equivalently written as a joint cumulative distribution function (CDF): 	 
\begin{equation}\label{eq1}
\begin{aligned}
\alpha_i &\leq \Pb\left[\, \bs{x}_i\in \mathbb{X}_i \ , \ \bs{x}_j\in  \tilde{\cl{B}}_j \, \right] \\
&=\int_{\mathbb{X}_i}\int_{\tilde{\cl{B}}_j} p(\bs{x}_i\,,\,\bs{x}_j) \, \dif \bs{x}_i \, \dif \bs{x}_j = F_{\bs{x}_i\,,\,\bs{x}_j}\,(\mathbb{X}_i \ , \ \tilde{\cl{B}}_j) \ ,
\end{aligned}
\end{equation}
where $p(\bs{x}_i\,,\,\bs{x}_j)$ is the joint probability density function (PDF) of $\bs{x}_i$ and $\bs{x}_j$.	 
Our goal is to calculate:
\begin{align*}
\Pb\left[\, \bs{x}_i\in \mathbb{X}_i \, \right] &=\int_{\mathbb{X}_i} p(\bs{x}_i) \, \dif \bs{x}_i  = F_{\bs{x}_i}\,(\mathbb{X}_i) \ ,
\end{align*}
where $p(\bs{x}_i)$ is the PDF of $\bs{x}_i$.
To transform the joint CDF into the marginal CDF of $\bs{x}_i$, one can take the limit of the joint CDF as $\tilde{\cl{B}}_j$ approaches $\R^{n_j}$:
\begin{equation}\label{eq2}
\begin{aligned}
\Pb\left[\, \bs{x}_i \in \mathbb{X}_i \, \right] &= F_{\bs{x}_i}\,(\mathbb{X}_i)  = \underset{\tilde{\cl{B}}_j \rightarrow \R^{n_j}}{\text{lim}} F_{\bs{x}_i\,,\,\bs{x}_j}\,(\mathbb{X}_i \ , \ \tilde{\cl{B}}_j) \\
&= \underset{\tilde{\cl{B}}_j \rightarrow \R^{n_j}}{\text{lim}} F_{\bs{x}_i\,|\,\bs{x}_j}\,(\mathbb{X}_i \ | \ \tilde{\cl{B}}_j) F_{\bs{x}_j}\,(\tilde{\cl{B}}_j) \\
&= F_{\bs{x}_i}\,(\mathbb{X}_i) \, \underset{\tilde{\cl{B}}_j \rightarrow \R^{n_j}}{\text{lim}} F_{\bs{x}_j}\,(\tilde{\cl{B}}_j) \ ,
\end{aligned}
\end{equation}
where the last equality is due to the fact that $\bs{x}_i$ and $\bs{x}_j$ are conditionally independent.

To determine $\underset{\tilde{\cl{B}}_j \rightarrow \R^{n_j}}{\text{lim}} F_{\bs{x}_j}\,(\tilde{\cl{B}}_j)$, one can calculate:
\vspace{-0.25cm}
\begin{equation}\label{eq3}
\begin{aligned}
\underset{\tilde{\cl{B}}_j \rightarrow \R^{n_j}}{\text{lim}} F_{\bs{x}_j}\,(\tilde{\cl{B}}_j) &= \int\limits_{\R^{n_j}} p(\bs{x}_j) \dif \bs{x}_j \\
&= \int\limits_{\R^{n_j} \setminus \tilde{\cl{B}}_j} p(\bs{x}_j) \dif \bs{x}_j + \int\limits_{\tilde{\cl{B}}_j} p(\bs{x}_j) \dif \bs{x}_j \\
&= \Pb\left[\, \bs{x}_j\notin \tilde{\cl{B}}_j \, \right] + \Pb\left[\, \bs{x}_j\in \tilde{\cl{B}}_j \, \right] \\
&= (1-\tilde{\alpha}_j) + \tilde{\alpha}_j = 1 \ ,
\end{aligned}
\end{equation}
where $p(\bs{x}_j)$ is the PDF of $\bs{x}_j$, and the last equality is a direct result of \Cref{reliable_thm}.
We now put all the steps together as follows:
\begin{align*}
\alpha_i&\leq \Pb\left[\, \bs{x}_i\in \mathbb{X}_i \ , \ \bs{x}_j\in  \tilde{\cl{B}}_j \, \right] = F_{\bs{x}_i\,,\,\bs{x}_j}\,(\mathbb{X}_i \ , \ \tilde{\cl{B}}_j)  \\
&\leq F_{\bs{x}_i}\,(\mathbb{X}_i) \, \underset{\tilde{\cl{B}}_j \rightarrow \R^{n_j}}{\text{lim}} F_{\bs{x}_j}\,(\tilde{\cl{B}}_j) \\
&= \Pb\left[\, \bs{x}_i \in \mathbb{X}_i \, \right] \left(\ \int\limits_{\R^{n_j} \setminus \tilde{\cl{B}}_j} p(\bs{x}_j) \dif \bs{x}_j + \int\limits_{\tilde{\cl{B}}_j} p(\bs{x}_j) \dif \bs{x}_j \ \right) \\
&\leq \int\limits_{\R^{n_j} \setminus \tilde{\cl{B}}_j} p(\bs{x}_j) \dif \bs{x}_j + \Pb\left[\, \bs{x}_i \in \mathbb{X}_i \, \right]  \int\limits_{\tilde{\cl{B}}_j} p(\bs{x}_j) \dif \bs{x}_j \\
&= (1-\tilde{\alpha}_j) + \tilde{\alpha}_j \ \Pb\left[\, \bs{x}_i \in \mathbb{X}_i \, \right]\ ,
\end{align*}
where the first inequality and equality is due to \eqref{eq1}, the second inequality is due to \eqref{eq2}, the second and last equality is due to \eqref{eq3}, and the last inequality is considering the worst-case situation, e.g. $\Pb \left[\, \bs{x}_i \in \mathbb{X}_i \, \big\vert \bs{x}_j \notin \tilde{\cl{B}}_j \,\right] = 1$.

By rearranging the last equation in above result: 
\begin{align}\label{res1_proof}
\frac{\alpha_i - (1-\tilde{\alpha}_j)}{\tilde{\alpha}_j} = 1 - \frac{1-\alpha_i}{\tilde{\alpha}_j} = \bar{\alpha}_i \ \leq \Pb\left[\, \bs{x}_i \in \mathbb{X}_i \, \right] \ .
\end{align}
This completes the proof of first part.
We now need to show the effect of having more than one neighboring agent.
To this end, the most straightforward step, in order to extend the current results, is to use the fact that all neighboring agents are independent from each other. 
We therefore can apply the previous results for a new situation where the agent $i$ with the probabilistic level of feasibility $\bar{\alpha}_i$ have another neighboring agent $\nu\in\cl{N}_i$ with $\tilde{\alpha}_\nu$ the level of reliability of $\tilde{\cl{B}}_\nu$.
By using \Cref{res1_proof}, we have the following relations for 	$j,\nu\in\cl{N}_i$
\begin{align*}
1 - \frac{1-\bar{\alpha}_i}{\tilde{\alpha}_\nu} 
&= 1 - \frac{1-\left(1 - \frac{1-\alpha_i}{\tilde{\alpha}_j}\right)}{\tilde{\alpha}_\nu} \\
&\qquad = 1 - \frac{1-{\alpha}_i}{\tilde{\alpha}_j \, \tilde{\alpha}_\nu} \ 
\leq \Pb\left[\, \bs{x}_i \in \mathbb{X}_i \, \right] \ .
\end{align*}
By continuing the similar arguments for all neighboring agents, one can obtain $\bar{\alpha}_i = 1 - \frac{1-\alpha_i}{\tilde{\alpha}_i} \leq \Pb\left[\, \bs{x}_i \in \mathbb{X}_i \, \right]$ such that $\tilde{\alpha}_i=  \tilde{\alpha}_1 \cdots \tilde{\alpha}_j \, \tilde{\alpha}_\nu \cdots \tilde{\alpha}_{|\cl{N}_i|}= \prod_{j\in\cl{N}_i}(\tilde{\alpha}_j)$.	
The proof is completed.
\QEDB

\bibliographystyle{IEEEtran}
\bibliography{dissertation}

\end{document}